%% file: run.tex
\documentclass[12pt,fleqn]{article}

\usepackage{amsmath,amsthm,amscd,amssymb}
\usepackage{latexsym}
\usepackage[all]{xy}

\allowdisplaybreaks

 \newcounter{abceqn}

 \newcounter{abcfig}

\newcommand{\re}{\operatorname{Re}}
\newcommand{\Span}{\operatorname{span}}
\newcommand{\norm}{\operatorname{norm}}

\newcommand{\vq}{\vec{q}}

\newcommand{\e}{\epsilon}

\newcommand{\W}{{\cal W}}
\renewcommand{\k}{\kappa}
\newcommand{\ga}{\gamma}
\newcommand{\Ga}{\Gamma}

\newcommand{\dl}{\delta}
\newcommand{\Dl}{\Delta}
\renewcommand{\th}{\theta}
\newcommand{\vth}{\vartheta}
\newcommand{\hvth}{\hat{\vartheta}}
\newcommand{\hy}{\hat{y}}
\newcommand{\htau}{\hat{\tau}}
\newcommand{\ra}{\rightarrow}
\newcommand{\al}{\alpha}
\newcommand{\be}{\beta}
\newcommand{\sg}{\sigma}

\newcommand{\pa}{\partial}

\newcommand{\la}{\lambda}
\newcommand{\bq}{\bar{q}}

\newcommand{\nid}{\noindent}

\newcommand{\om}{\omega}

\renewcommand{\theequation}{\thesection.\arabic{equation}}
\newcommand{\eqnsection}[1]{
	\section{#1}
	\setcounter{equation}{0}
	\renewcommand{\theequation}{\thesection.\arabic{equation}}
	\setcounter{figure}{0}
	\renewcommand{\thefigure}{\thesection.\arabic{figure}}
	\setcounter{remark}{0}
	\renewcommand{\theremark}{\thesection.\arabic{remark}}
	\setcounter{theorem}{0}
	\renewcommand{\thetheorem}{\thesection.\arabic{theorem}}
	\setcounter{lemma}{0}
	\renewcommand{\thelemma}{\thesection.\arabic{lemma}}
	\setcounter{proposition}{0}
	\renewcommand{\theproposition}{\thesection.\arabic{proposition}}
}

\title{\bf Persistent Homoclinic Orbits for Nonlinear Schr\"odinger 
Equation Under Singular Perturbation}

\author{\\ \\ \\ \\ Yanguang\ (Charles) \ \ Li  
\thanks{This work is partially supported by a Guggenheim Fellowship.}
\\ \\ \\ Department of Mathematics \\ \\ University of Missouri - Columbia
\\ \\ Columbia, MO 65211}

\date{\today}

\renewcommand{\theequation}{\thesection.\arabic{equation}}

\begin{document}
\bibliographystyle{plain}
\maketitle

\newpage
\begin{abstract} 
Existence of homoclinic orbits in the cubic nonlinear Schr\"odinger equation 
under singular perturbations is proved. Emphasis is placed upon the 
regularity of the semigroup $e^{\e t \pa_x^2}$ at $\e = 0$. This article is 
a substantial generalization of \cite{LMSW96}, and motivated by the effort 
of Dr. Zeng \cite{Zen00a} \cite{Zen00b}. The mistake of Zeng in \cite{Zen00b} 
is corrected with a normal form transform approach. Both one and two unstable 
modes cases are investigated.  
\end{abstract}

\newtheorem{lemma}{Lemma}
\newtheorem{theorem}{Theorem}
\newtheorem{corollary}{Corollary}
\newtheorem{remark}{Remark}
\newtheorem{definition}{Definition}
\newtheorem{proposition}{Proposition}
\newtheorem{assumption}{Assumption}

\newpage
\tableofcontents

\newpage
\input intro

\newpage
\input local

\newpage
\input global

\newpage
\bibliography{NLS}

\end{document}

%% file: intro.tex
\eqnsection{Introduction}

Consider the singularly perturbed nonlinear Schr\"odinger equation,
\begin{equation}
iq_t = q_{xx} +2 [ |q|^2 - \om^2] q +i \e [q_{xx} - \al q +\be ] \ ,
\label{pnls}
\end{equation}
where $q = q(t,x)$ is a complex-valued function of the two real 
variables $t$ and $x$, $t$ represents time, and $x$ represents
space. $q(t,x)$ is subject to periodic boundary condition of period 
$2 \pi$, and even constraint, i.e. 
\[
q(t,x + 2 \pi) = q(t,x)\ , \ \ q(t,-x) = q(t,x)\ .
\]
$\om$ is a positive constant, $\al >0$ and $\be >0$ are constants, 
and $\e > 0$ is the perturbation parameter.

In this work, we revisit the problem on the existence of homoclinic 
orbits in perturbed nonlinear Schr\"odinger equations \cite{LMSW96}.
The crucial new feature is as follows:
Singular perturbation $\e \pa_x^2 q$ will be investigated in 
contrast to the regular perturbation $\e \hat{\pa}_x^2 q$ considered 
in \cite{LMSW96}, where $\hat{\pa}_x^2$ is a pseudo-differential 
operator obtained from a truncation of $\pa_x^2$. This study is motivated 
by the effort of Dr. Zeng \cite{Zen00a} \cite{Zen00b}.
The main difficulty introduced by the singular 
perturbation $\e \pa_x^2$ is that it breaks the spectral gap condition 
of the unperturbed system. Therefore, standard invariant manifold results 
will not apply. Nevertheless, it turns out that certain invariant manifold 
results do hold. The regularity of such invariant manifolds at $\e = 0$ is 
controled by the regularity of $e^{\e \pa_x^2}$ at $\e = 0$. Difficulties 
and interesting results created by the singular perturbation term 
$\e \pa_x^2 q$ will all be commented in Remarks.

The entire theory of locating a homoclinic orbit is divided into two parts.
Part 1 deals with local invariant manifold theory. Part 2 deals with global 
theory which includes integrable theory, Melnikov analysis, etc..

The notation $| \ |$ will denote absolute value, and the notation $\| \ \|_s$ 
will denote the Sobolev $H^s$ (i.e. $W^{s,2}$) norm of periodic function 
with period $2\pi$. 

This article is written for experts only. Standard details will be omitted.

%% file: local.tex
\eqnsection{Local Theory}

Local theory is referred to a theory in a neighborhood of certain 
circle of fixed points, which includes local unstable fiber theorem,
local center-stable manifold theorem, and size estimate of local 
stable manifold for certain saddle. These are some of the tools needed
in locating a homoclinic orbit.

\subsection{Dynamics in a 2D Invariant Subspace}

The 2D subspace $\Pi$,
\begin{equation}\label{Pi} \Pi=\{ q\mid \ \partial_xq=0\},\end{equation}
is an invariant subspace under the PNLS flow (\ref{pnls}). The 
governing equation in $\Pi$ is
\begin{equation} i\dot{q}=2[|q|^2-\om^2]q+i\epsilon [-\alpha 
q+\beta],\label{Pie1}\end{equation}
where $\cdot =\frac{d}{dt}\ $.
Dynamics of this equation is shown in Figure~\ref{figPi}. 
Interesting dynamics is created through resonance in the neighborhood of the
circle $S_\om$:
\begin{equation}\label{rcl} S_\om=\{ q\in \Pi \mid \ |q|=\om \}.\end{equation}
When $\epsilon =0$, $S_\om$ consists of fixed points. To explore the 
dynamics in this neighborhood better, one can make a series of
changes of coordinates. Let $q=\sqrt{I}e^{i\theta}$,
then \eqref{Pie1} can be rewritten as
\begin{align} \dot{I}&= \epsilon ( -2\alpha I+2\beta \sqrt{I}\cos 
\theta )\ ,\label{Ithe1}\\
\dot{\theta} &=-2(I-\om^2)-\epsilon \beta \frac{\sin 
\theta}{\sqrt{I}}\ .\label{Ithe2}\end{align}
There are three fixed points:
\begin{enumerate}\item The focus $O_\epsilon$ in the neighborhood of 
the origin,
\begin{equation}\begin{cases} I=\epsilon ^2\frac{\beta^2}{4\omega^4}+\cdots ,\\
\cos \theta=\frac{\alpha \sqrt{I}}{\beta}, & \theta \in \left( 
0,\frac{\pi}{2}\right).\end{cases}\label{Oec}\end{equation}
Its eigenvalues are
\begin{equation}\label{Oee} \mu_{1,2}=\pm 
i\sqrt{4(\omega^2-I)^2-4\epsilon \sqrt{I}\beta \sin \theta}-\epsilon 
\alpha,\end{equation}
where $I$ and $\theta$ are given in \eqref{Oec}.
\item The focus $P_\epsilon$ in the neighborhood of $S_\omega$ \eqref{rcl},
\begin{equation}\label{Pec}\begin{cases} I=\omega^2+\epsilon 
\frac{1}{2\omega}\sqrt{\beta^2-\alpha^2\omega^2}+\cdots ,\\
\cos \theta=\frac{\alpha \sqrt{I}}{\beta},& \theta \in \left( 
-\frac{\pi}{2}, 0\right).\end{cases}\end{equation}
Its eigenvalues are
\begin{equation}\label{Pee}\mu_{1,2}=\pm 
i\sqrt{\epsilon}\sqrt{-4\sqrt{I}\beta \sin \theta+\epsilon \left( 
\frac{\beta \sin
\theta}{\sqrt{I}}\right)^2}-\epsilon \alpha,\end{equation}
where $I$ and $\theta$ are given in \eqref{Pec}.
\item The saddle $Q_\epsilon$ in the neighborhood of $S_\omega$ \eqref{rcl},
\begin{equation}\label{Qec}\begin{cases}I=\omega^2-\epsilon 
\frac{1}{2\omega}\sqrt{\beta^2-\alpha^2\omega^2}+\cdots ,\\
\cos \theta  =\frac{\alpha \sqrt{I}}{\beta},&\theta \in \left( 
0,\frac{\pi}{2}\right).\end{cases}\end{equation}
Its eigenvalues are
\begin{equation}\label{Qee} \mu_{1,2}=\pm 
\sqrt{\epsilon}\sqrt{4\sqrt{I}\beta\sin \theta -\epsilon \left( 
\frac{\beta \sin
\theta}{\sqrt{I}}\right)^2}-\epsilon \alpha,\end{equation}
where $I$ and $\theta $ are given in \eqref{Qec}.
\end{enumerate}
\begin{figure}
\label{figPi}
\vspace{1.5in}
\caption{Dynamics on the invariant plane $\Pi$.}
\end{figure}
\nid
Now focus our attention to order $\sqrt{\epsilon}$ neighborhood of 
$S_\omega$ \eqref{rcl} and let
\begin{equation*}J=I-\omega^2,\quad J=\sqrt{\epsilon}j,\quad 
\tau=\sqrt{\epsilon}t,\end{equation*}
we have
\begin{align}\label{je1} j'&= 2\left[ -\alpha 
(\omega^2+\sqrt{\epsilon}j) +\beta 
\sqrt{\omega^2+\sqrt{\epsilon}j}\cos \theta \right],\\
\label{je2}\theta '&= -2j-\sqrt{\epsilon}\beta \frac{\sin 
\theta}{\sqrt{\omega^2+\sqrt{\epsilon}j}},\end{align}
where $'=\frac{d}{d\tau}\ $. To leading order, we get
\begin{align}\label{je3} j'&=2[-\alpha \omega^2+\beta \omega \cos \theta]\ ,\\
\label{je4}\theta' & =-2j\ .\end{align}
There are two fixed points which are the counterparts of $P_\epsilon$ 
and $Q_\epsilon$ \eqref{Pec} and \eqref{Qec}:
\begin{enumerate}\item The center $P_*$,
\begin{equation}\label{Pc} j=0,\quad \cos \theta=\frac{\alpha 
\omega}{\beta},\quad \theta \in
\left(-\frac{\pi}{2},0\right).\end{equation}
Its eigenvalues are
\begin{equation}\label{Pe} \mu_{1,2}=\pm 
i2\sqrt{\omega}(\beta^2-\alpha^2\omega^2)^{\frac{1}{4}}.\end{equation}
\item The saddle $Q_*$,
\begin{equation}\label{Qc} j=0,\quad \cos \theta=\frac{\alpha 
\omega}{\beta},\quad \theta \in \left(
0,\frac{\pi}{2}\right).\end{equation} Its eigenvalues are
\begin{equation}\label{Qe} \mu_{1,2}=\pm 
2\sqrt{\omega}(\beta^2-\alpha ^2\omega^2)^{\frac{1}{4}}.\end{equation}
\end{enumerate}
\nid
In fact, \eqref{je3} and \eqref{je4} form a Hamiltonian system with the 
Hamiltonian
\begin{equation}\label{fham} \mathcal{H}=j^2+2\omega (-\alpha \omega 
\theta+\beta \sin \theta).\end{equation}
Connecting to $Q_*$ is a fish-like singular level set of 
$\mathcal{H}$, which intersects the axis $j=0$ at $Q_*$ and
$\hat{Q}=(0,\hat{\theta})$,
\begin{equation}\label{head} \alpha \omega 
(\hat{\theta}-\theta_*)=\beta (\sin \hat{\theta}-\sin \theta_*),\quad 
\hat{\theta}\in (-\frac{3\pi}{2}, 0),\end{equation}
where $\theta_*$ is given in \eqref{Qc}. See Figure~\ref{fish} for an 
illustration of the dynamics of \eqref{je1}-\eqref{je4}. 
\begin{figure}
\vspace{1.5in}
\caption{The fish-like dynamics in the neighborhood of the resonant circle
$S_\om$.}
\label{fish}
\end{figure}
\nid
For later use, we define
a piece of each of the stable and unstable manifolds of $Q_*$,
\begin{equation*}j=\phi_*^u(\theta),\quad j=\phi^s_*(\theta),\quad 
\theta\in [\hat{\theta}+\hat{\delta}, \theta_*+2\pi],\end{equation*}
for some small $\hat{\delta}>0$, and
\begin{equation}\begin{split}\label{cur} 
\phi^u_*(\theta)&=-\frac{\theta-\theta_*}{|\theta 
-\theta_*|}\sqrt{2\omega [\alpha \omega
(\theta-\theta_*)-\beta (\sin \theta-\sin \theta_*)]},\\
 \phi^s_*(\theta )&=-\phi^u_*(\theta).\end{split}\end{equation}
$\phi^u_*(\theta)$ and $\phi^s_*(\theta)$ perturb smoothly in $\theta 
$ and $\sqrt{\epsilon}$ into $\phi^u_{\sqrt{\epsilon}}$ and 
$\phi^s_{\sqrt{\epsilon}}$ for
\eqref{je1} and \eqref{je2}.

The homoclinic orbit to be located will take off from $Q_\epsilon$ 
along its unstable curve, flies away from and returns to $\Pi$, lands
near the stable curve of $Q_\epsilon$ and approaches $Q_\epsilon$ spirally.

\subsection{Change of Coordinates}

As mentioned above, interesting dynamics happens in the neighborhood 
of the circle $S_\omega$ \eqref{rcl}. It is natural and convenient
to center our coordinates around $S_\omega$. First, write $q$ as
\begin{equation}\label{pcd} q(t,x)=[\rho(t)+f(t,x)]e^{i\theta (t)},\end{equation}
where $\rho$ and $\theta$ are polar coordinates on $\Pi $ \eqref{Pi}, 
and $f$ has zero spatial mean. We use the notation $\langle \cdot
\rangle$ to denote spatial mean,
\begin{equation}\label{mean} \langle q\rangle 
=\frac{1}{2\pi}\int^{2\pi}_0qdx.\end{equation}
Since the $L^2$-norm is an action variable when $\epsilon =0$, it is 
more convenient to replace $\rho$ by:
\begin{equation}\label{L2n} I=\langle |q|^2\rangle =\rho^2+\langle 
|f|^2\rangle.\end{equation}
Since $S_\om$ corresponds to $I=\omega^2$, the final pick is
\begin{equation}\label{Jpc}J=I-\omega^2.\end{equation}
In terms of the new variables $(J,\theta, f)$, Equation \eqref{pnls} 
can be rewritten as
\begin{align}\label{nc1} \dot{J}&=\epsilon \left[-2\alpha 
(J+\omega^2)+2\beta \sqrt{J+\omega^2}\cos \theta \right]+\epsilon
\mathcal{R}^J_2,\\
\label{nc2} \dot{\theta}&=-2J-\epsilon \beta \frac{\sin 
\theta}{\sqrt{J+\omega^2}}+\mathcal{R}^\theta_2,\\
\label{nc3}f_t&=L_\epsilon f+V_\epsilon 
f-i\mathcal{N}_2-i\mathcal{N}_3,\end{align}
where
\begin{align}\label{wnc1} L_\epsilon f&=-if_{xx}+\epsilon (-\alpha 
f+f_{xx})-i2\omega ^2(f+\bar{f}),\\
\label{wnc2}V_\epsilon f&=-i2J(f+\bar{f})+i\epsilon \beta f\frac{\sin 
\theta}{\sqrt{J+\omega^2}},\\
\label{wnc3}\mathcal{R}^J_2&=-2\langle |f_x|^2\rangle +2\beta \cos 
\theta \left[ \sqrt{J+\omega^2-\langle
|f|^2\rangle}-\sqrt{J+\omega^2}\right],\\
\begin{split}\label{wnc4}\mathcal{R}^\theta_2&=-\langle 
(f+\bar{f})^2\rangle -\frac{1}{\rho}\langle |f|^2(f+\bar{f})\rangle\\
&\quad  -\epsilon
\beta\sin
\theta
\left[
\frac{1}{\sqrt{J+\omega^2-\langle 
|f|^2\rangle}}-\frac{1}{\sqrt{J+\omega^2}}\right],\end{split}\\
\label{wnc5} \mathcal{N}_2&=2\rho [2(|f|^2-\langle |f|^2\rangle 
)+(f^2-\langle f^2\rangle )],\\
\begin{split}\label{wnc6} \mathcal{N}_3&=-\langle 
f^2+\bar{f}^2+6|f|^2\rangle f+2(|f|^2f-\langle |f|^2f\rangle )\\
&\quad -\frac{1}{\rho} \langle |f|^2(f+\bar{f})\rangle f-2\langle 
|f|^2\rangle \bar{f}\\
&\quad -\epsilon \beta \sin \theta \left[ \frac{1}{\sqrt{J+\omega^2-\langle
|f|^2\rangle}}-\frac{1}{\sqrt{J+\omega^2}}\right]f.\end{split}\end{align}

\begin{remark} The singular perturbation term ``$\epsilon \pa_x^2q$" 
can be seen at two locations, $L_\epsilon$ and $\mathcal{R}^J_2$
(\ref{wnc1},\ref{wnc3}). The singular perturbation term $\langle 
|f_x|^2\rangle $ in $\mathcal{R}^J_2$ does not create any difficulty.
Since $H^1$ is a Banach algebra~\cite{Ada75}, this term is still of 
quadratic order, $\langle |f_x|^2\rangle \sim \mathcal{O}(\|
f\|^2_1)$.\end{remark}
\begin{lemma} The nonlinear terms have the orders:
\begin{equation*}\begin{split}&|\mathcal{R}^J_2|\sim \mathcal{O}(\| 
f\|^2_s),\quad |\mathcal{R}^\theta_2|\sim \mathcal{O}(\|
f\|^2_s),\\
& \|\mathcal{N}_2\|_s\sim
\mathcal{O} (\| f\|^2_s),\quad \| \mathcal{N}_3\|_s\sim 
\mathcal{O}(\| f\|^3_s),\quad (s\geq 
1).\end{split}\end{equation*}\end{lemma}

\begin{proof} The proof is an easy direct verification.\end{proof}

\subsection{Normal Form Transformation}

In locating a homoclinic orbit to $Q_\epsilon$ \eqref{Qec}, we need 
to estimate the size of the local stable manifold of $Q_\epsilon$.
The size of the variable $J$ is of order 
$\mathcal{O}(\sqrt{\epsilon})$. The size of the variable $\theta$ is 
of order $\mathcal{O}(1)$.
To be able to track a homoclinic orbit, we need the size of the 
variable $f$ to be of order $\mathcal{O}(\epsilon^\mu)$, $\mu <1$. 
Such
an estimate can be achieved, if the quadratic term 
$\mathcal{N}_2$ \eqref{wnc5} in \eqref{nc3} can be removed through a 
normal form
transformation.

In terms of Fourier transforms,
\begin{equation*}f=\sum_{k\neq0}\hat{f}(k)e^{ikx},\quad 
\bar{f}=\sum_{k\neq 0}\overline{\hat{f}(-k)}e^{ikx},\end{equation*}
and the two terms in $\mathcal{N}_2$ can be written as,
\begin{align} \notag f^2-\langle f^2\rangle &= \sum_{k+\ell \neq 
0}\hat{f}(k)\hat{f}(\ell)e^{i(k+\ell)x},\\
\begin{split}\label{sytr} |f|^2-\langle |f|^2\rangle &= 
\sum_{k+\ell\neq 0}\hat{f}(k)\overline{\hat{f}(-\ell)}e^{i(k+\ell)x}\\
&=\frac{1}{2}\sum_{k+\ell \neq
0}[\hat{f}(k)\overline{\hat{f}(-\ell)}+\hat{f}(\ell)\overline{\hat{f}(-k)} 
]e^{i(k+\ell)x}.\end{split}\end{align}
It turns out to be convenient to work with the symmetrized form 
\eqref{sytr}. We will search for a normal form transformation of the 
general form,
\begin{equation}\label{nft} g=f+K(f,f),\end{equation}
where
\begin{equation*}\begin{split} K(f,f)&= \sum_{k+\ell\neq 0}\left[
\hat{K}_1(k,\ell)\hat{f}(k) 
\hat{f}(\ell)+\hat{K}_2(k,\ell)\hat{f}(k)\overline{\hat{f}(-\ell)}\right.\\
&\quad \left.+\hat{K}_2(\ell,k)\overline{\hat{f}(-k)}\hat{f}(\ell)+
\hat{K}_3(k,\ell)\overline{\hat{f}(-k)}\overline{\hat{f}(-\ell)}\right]e^{i(k+\ell)x},\end{split}\end{equation*}
$\hat{K}_j(k,\ell)$, $(j=1,2,3)$ are the unknown coefficients to 
be determined, and 
$\hat{K}_j(k,\ell)=\hat{K}_j(\ell,k)$, $(j=1,3)$.

\begin{lemma} For $\omega \in 
\left(\frac{1}{2},\frac{3}{2}\right)/S$, $S$ is a finite subset, 
there exists a normal form transformation
of the form \eqref{nft} that transforms the equation
\begin{equation*}f_t=L_\epsilon f-i\tilde{\mathcal{N}}_2,\end{equation*}
into an equation with a cubic nonlinearity
\begin{equation*}g_t=L_\epsilon g+\mathcal{O}(\| g\|^3_s),\quad 
(s\geq 1),\end{equation*}
where $L_\epsilon$ is given in \eqref{wnc1}, and 
$\tilde{\mathcal{N}}_2$ has the expression (cf: \eqref{wnc5}),
\begin{equation}\label{nnf} \tilde{\mathcal{N}}_2=2\omega 
[2(|f|^2-\langle |f|^2\rangle )+(f^2-\langle f^2\rangle 
)].\end{equation}
\end{lemma}

\begin{proof} Denote the operator $i\partial _t-iL_\epsilon$ by 
$\mathcal{L}_\epsilon$. We have
\begin{equation*}\begin{split}\mathcal{L}_\epsilon 
g&=\mathcal{L}_\epsilon f+\mathcal{L}_\epsilon K(f,f)\\
&=\mathcal{L}_\epsilon f-iL_\epsilon K(f,f)+iK(L_\epsilon 
f,f)+iK(f,L_\epsilon f)\\
&+iK(\partial_tf-L_\epsilon f,f)+iK(f,\partial_tf-L_\epsilon 
f),\end{split}\end{equation*}
where $\mathcal{L}_\epsilon f=\tilde{\mathcal{N}}_2$ and 
$K(\partial_tf-L_\epsilon f,f)$ and $K(f,\partial_tf-L_\epsilon f)$ 
will be
shown to be cubic in $f$. To eliminate the quadratic terms, we need to set
\begin{equation*}iL_\epsilon K(f,f)-iK(L_\epsilon 
f,f)-iK(f,L_\epsilon f)=\tilde{\mathcal{N}}_2,\end{equation*}
which takes the explicit form:
\begin{eqnarray}
& &(\sigma_1+i\sigma)\hat{K}_1(k,\ell)+b\hat{K}_2(k,\ell)+b\hat{K}_2(\ell,k)+b\overline{\hat{K}_3(-k,-\ell)}=-2\omega,\label{nore1}\\
& &-b\hat{K}_1(k,\ell)+(\sigma_2+i\sigma)\hat{K}_2(k,\ell)+b 
\overline{\hat{K}_2(-\ell,-k)}+b\hat{K}_3(k,\ell)=-2\omega
,\label{nore2}\\
& &-b\hat{K}_1(k,\ell)+b
\overline{\hat{K}_2(-k,-\ell)}+(\sigma_3+i\sigma)\hat{K}_2(\ell,k)+b\hat{K}_3(k,\ell)= 
-2\omega
,\label{nore3}\\
& &b \overline{\hat{K}_1(-k,-\ell)}-b\hat{K}_2(k,\ell) 
-b\hat{K}_2(\ell,k)+(\sigma _4+i\sigma
)\hat{K}_3(k,\ell)=0,\label{nore4}
\end{eqnarray}
where
\begin{align*} b&=-2\omega^2,\quad \sigma=\epsilon 
(2k\ell -\alpha),\quad \sigma_1=2(k\ell+\omega^2),\quad
\sigma_2=2(\ell^2+k\ell-\omega^2),\\
\sigma_3&=2(k^2+k\ell-\omega^2),\quad 
\sigma_4=2(k^2+\ell^2+k\ell-3\omega^2).\end{align*}
Since these coefficients are even in $(k,\ell)$, we will search for 
even solutions, i.e.
\begin{equation*}\hat{K}_j(-k,-\ell)=\hat{K}_j(k,\ell),\quad 
j=1,2,3.\end{equation*}
\eqref{nore1}+\eqref{nore4}, \eqref{nore2}+\eqref{nore4}$^-$, and 
\eqref{nore3}+\eqref{nore4}$^-$ lead to
\begin{align*} 
(\sigma_1+i\sigma)\hat{K}_1(k,\ell)+b\overline{\hat{K}_1(k,\ell)}&=-K,\\
(\sigma_2+i\sigma)\hat{K}_2(k,\ell)-b 
\overline{\hat{K}_2(k,\ell)}&=-\bar{K},\\
(\sigma_3+i\sigma)\hat{K}_2(\ell, 
k)-b\overline{\hat{K}_2(\ell,k)}&=-\bar{K},\end{align*}
where
\begin{equation*}K=2\omega+(\sigma_4+i\sigma)\hat{K}_3(k,\ell)+b\overline{\hat{K}_3(k,\ell)}.\end{equation*}
Therefore we can express $\hat{K}_j(k,\ell)$ in terms of $K$ as,
\begin{align}\label{nors1} \hat{K}_1(k,\ell)&= 
(\sigma^2_1+\sigma^2-b^2)^{-1}[b\bar{K}-(\sigma_1-i\sigma)K],\\
\label{nors2}\hat{K}_2(k,\ell)&=(\sigma^2_2+\sigma^2-b^2)^{-1}[-b\bar{K}-(\sigma_2-i\sigma)K],\\
\label{nors3} \hat{K}_2(\ell,k)&= 
(\sigma^2_3+\sigma^2-b^2)^{-1}[-b\bar{K}-(\sigma_3-i\sigma)K],\\
\label{nors4} 
\hat{K}_3(k,\ell)&=(\sigma^2_4+\sigma^2-b^2)^{-1}[(\sigma_4-i\sigma)(K-2\omega)-b(\bar{K}-2\omega)].\end{align}
Substituting these expressions into \eqref{nore4}, we get
\begin{equation}\label{nors5} 
K=(|U|^2-|V|^2)^{-1}(W\bar{U}-\bar{W}V),\end{equation}
where
\begin{align}\begin{split}\label{nors6}
U&=\frac{b^2}{\sigma^2_1+\sigma^2-b^2}+\frac{b(\sigma_2-i\sigma)}{\sigma^2_2+\sigma^2-b^2}+
\frac{b(\sigma_3-i\sigma)}{\sigma^2_3+\sigma^2-b^2}
+\frac{\sigma^2_4+\sigma^2}{\sigma^2_4+\sigma^2-b^2},\end{split}\\
\begin{split}\label{nors7}
V&=-\frac{b(\sigma_1+i\sigma)}{\sigma^2_1+\sigma^2-b^2}
+\frac{b^2}{\sigma^2_2+\sigma^2-b^2}+
\frac{b^2}{\sigma^2_3+\sigma^2-b^2}-\frac{b(\sigma_4+i\sigma)}{\sigma^2_4+\sigma^2-b^2},\end{split}\\
\label{nors8} W&=2\omega 
(\sigma^2_4+\sigma^2-b^2)^{-1}[\sigma^2_4+\sigma^2-b(\sigma_4+i\sigma)].\end{align}
For $\omega\in \left(\frac{1}{2},\frac{3}{2}\right)$, the 
denominators in \eqref{nors1}-\eqref{nors4} and 
\eqref{nors6}-\eqref{nors8},
and $\sigma_j \ (1\leq j\leq 4)$ vanish at $\omega$ in a finite subset. 
For $\omega$ not in this finite subset, $\sigma^2$ is a
$\mathcal{O}(\epsilon^2)$ small perturbation of $\sigma_j^2-b^2 \ (1\leq 
j\leq 4)$, and $\sigma$ is a $\mathcal{O}(\epsilon)$ small
perturbation of $\sigma_j \ (1\leq j\leq 4)$. Setting $\sigma=0$, we 
have $K=2\omega$, which leads to the solution given in \cite{LMSW96}.
To the order $\mathcal{O}(\epsilon)$,
\begin{equation}\begin{split}\label{nors9} K&=2\omega \left[ 1+ib\sigma \left(
\frac{1}{\sigma^2_1-b^2}+\frac{1}{\sigma^2_2-b^2}+\frac{1}{\sigma^2_3-b^2}\right) 
\cdot \right.\\
&\quad \cdot \left. \left( \frac{b}{\sigma_1-b}+\frac{b}{\sigma_2+b}
+\frac{b}{\sigma_3+b}+\frac{\sigma_4}{\sigma_4-b}\right)^{-1}\right].\end{split}\end{equation}
We will show that the denominator in \eqref{nors9} does not vanish 
except for $\omega$ in a finite subset. Denote the denominator by
$D$. As $k\to \pm \infty$, or $\ell \to \pm \infty$,
\begin{equation}\label{nors10} D\to 1.\end{equation}
We also know that
\begin{equation*}\begin{split} D&=1+b\left[ 
\frac{1}{\sigma_1-b}+\frac{1}{\sigma_2+b}+\frac{1}{\sigma_3+b}+\frac{1}{\sigma_4-b}\right]\\
&=1+\frac{b}{2} \left[ 
\frac{(k+\ell)^2-4\omega^2}{(\ell^2+k\ell-2\omega^2)(k^2+k\ell 
-2\omega^2)}
+\frac{(k+\ell)^2}{(k\ell+2\omega^2)(k^2+\ell^2+k\ell-2\omega^2)}\right]\\
&=\frac{-2(2\omega^2)^4+\chi
_1(2\omega^2)^3+\chi_2(2\omega^2)^2+\chi_3(2\omega^2)+\chi_4}{(\ell^2+k\ell-2\omega^2)
(k^2+k\ell-2\omega^2)(k\ell+2\omega^2)(k^2+\ell^2+k\ell-2\omega^2)},\end{split}\end{equation*}
where $\chi_j \ (1\leq j\leq 4)$ are polynomials in $k$ and $\ell$. For 
each $k$ and $\ell$, the numerator vanishes at most at four values
of $2\omega^2$. Together with the fact \eqref{nors10}, we have that 
for $\omega \in \left(\frac{1}{2},\frac{3}{2}\right)$, $D$ does not
vanish except for $\omega$ in a finite subset.

The denominator in \eqref{nors5} has the representation:
\begin{equation*}|U|^2-|V|^2=\re \{ (U+V)(\bar{U}-\bar{V})\},\end{equation*}
where
\begin{align*} U+V&=\frac{\sigma_4}{\sigma_4+b}-ib\sigma \left[
\frac{1}{\sigma^2_1-b^2}+\frac{1}{\sigma^2_2-b^2}+\frac{1}{\sigma^2_3-b^2}+\frac{1}{\sigma^2_4-b^2} 
\right]\\
&\quad + \text{higher order terms in }\epsilon,\\
U-V&=D+ib\sigma \left[ 
\frac{1}{\sigma^2_1-b^2}-\frac{1}{\sigma^2_2-b^2}-\frac{1}{\sigma^2_3-b^2}+\frac{1}{\sigma^2_4-b^2}\right]\\
&\quad + \text{ higher order terms in }\epsilon.\end{align*}
Then
\begin{equation*}|U|^2-|V|^2=\frac{\sigma_4}{\sigma_4+b}D+\mathcal{O}(\epsilon^2).\end{equation*}
Therefore, for $\omega \in \left( \frac{1}{2},\frac{3}{2}\right)$, 
the denominator in \eqref{nors5}, $|U|^2-|V|^2$ does not vanish
except for $\omega$ in a finite subset. \eqref{nors1}-\eqref{nors5} 
give the solution to the linear system \eqref{nore1}-\eqref{nore4}
for $\omega \in \left( \frac{1}{2},\frac{3}{2}\right)/S$, where $S$ 
is a finite subset.

As in \cite{LMSW96}, $K(f,f)$ is also a bounded bilinear map:
\begin{equation*}\| K(f,f)\|_s\leq C\| f\|^2_s,\quad (s\geq 1).\end{equation*}
We can invert the equation \begin{equation*}g=f+K(f,f)\end{equation*}
to obtain
\begin{equation*}f=g+\mathcal{K}(g),\end{equation*}
where $\mathcal{K}$ is of order $\mathcal{O}(\| g\|^2_s)$, $(s\geq 
1)$. Thus, terms like $K(\pa_tf-L_\epsilon f,f)$ and
$K(f,\pa_tf-L_\epsilon f)$ are cubic terms in $g$.
\end{proof}

\begin{remark} In this remark, we would like to make a comparison 
between the above normal form transform with that in \cite{LMSW96},
and in particular to comment on why the above normal form transform 
is necessary when singular perturbation $\epsilon \partial^2_x f$
is studied. In \cite{LMSW96}, the linear operator $L_\epsilon$ is 
replaced by $L_0$ (i.e. setting $\epsilon=0$ in $L_\epsilon$) in
constructing normal form transform. The corresponding normal form transform 
is given by
\begin{equation}\label{nft0} K=2\omega,\quad 
\hat{K}_3(k,\ell)=0,\quad 
\hat{K}_1(k,\ell)=-\frac{\omega}{k\ell},\quad
\hat{K}_2(k,\ell)=-\frac{\omega}{\ell(k+\ell)}.\end{equation}
When such a normal form transform is applied to Equation \eqref{nc3}, 
the singular perturbation term $\epsilon \partial^2_x f$ will
introduce the following term in the equation for $g$:
\begin{equation}\label{indu} \epsilon \partial^2_x 
\mathcal{K}(f,f)\end{equation}
which is actually an unbounded bilinear operator. Therefore, we have to 
work with 
$L_\epsilon$ for a normal form transform.
On the other hand, in \cite{LMSW96}, the singular perturbation 
$\epsilon \partial^2_x$ is mollified into a bounded
pseudo-differential operator (actually a bounded Fourier multiplier) 
$\epsilon\hat{\partial}^2_x$. The term \eqref{indu} is replaced
by
\begin{equation*}\epsilon \hat{\partial}^2_x \mathcal{K}(f,f)\end{equation*}
which is of order $\mathcal{O}(\epsilon \| f\|_s^2)$ sufficient for the 
estimate on the size of the local stable manifold of $Q_\epsilon$.
Although the normal form transform \eqref{nors1}-\eqref{nors5} has a 
more complicated expression than \eqref{nft0}, they have the same
asymptotic nature in $k$ and $\ell$. $\epsilon \partial^2_x f$ 
only introduces small perturbations in the expressions of
$\hat{K}_j(k,\ell) \ (1\leq j\leq 3)$.\end{remark}

We apply the normal form transform given by 
\eqref{nors1}-\eqref{nors5} to the full equation \eqref{nc3}, and the 
full system
\eqref{nc1}-\eqref{nc3} is transformed into:
\begin{align}\label{nfe1} \dot{J}&=\epsilon \left[ -2\alpha 
(J+\omega^2)+2\beta \sqrt{J+\omega^2}\cos \theta \right] +\epsilon
\mathcal{R}^J_2,\\
\label{nfe2} \dot{\theta} &=-2J-\epsilon\beta \frac{\sin 
\theta}{\sqrt{J+\omega^2}}+\mathcal{R}^\theta_2,\\
\label{nfe3} g_t&=L_\epsilon g+V_\epsilon g+\mathcal{N},\end{align}
where $L_\epsilon$, $V_\epsilon$, $\mathcal{R}^J_2$ and 
$\mathcal{R}^\theta_2$ are given in \eqref{wnc1}-\eqref{wnc4} with
$f=g+\mathcal{K}(g)$, and
\begin{equation}\begin{split}\label{wnf} \mathcal{N}&=V_\epsilon \mathcal{K}
(g)-i(\mathcal{N}_2-\tilde{\mathcal{N}}_2)-i\mathcal{N}_3+K(\partial 
_tf-L_\epsilon f,f)+ \mathcal{K}(f,\partial_tf-L_\epsilon f)\\
&=V_\epsilon 
\mathcal{K}(g)-i(\mathcal{N}_2-\tilde{\mathcal{N}_2})-i\mathcal{N}_3
+K(V_\epsilon f-i\mathcal{N}_2-i\mathcal{N}_3,f)\\
& \quad +K(f,V_\epsilon
f-i\mathcal{N}_2-i\mathcal{N}_3),\end{split}\end{equation}
where $\mathcal{N}_2$, $\tilde{\mathcal{N}}_2$ and $\mathcal{N}_3$ 
are given in \eqref{wnc5}, \eqref{nnf} and \eqref{wnc6} with
$f=g+\mathcal{K}(g)$. $\mathcal{N}$ has the estimate,
\begin{equation}\label{wnfe} \| \mathcal{N}\|_s\sim \mathcal{O} (|J 
|\|g\|^2_s+\epsilon \| g\|^2_s+\| g\|^3_s),\quad
(s\geq 1).\end{equation}

\subsection{Unstable Fibers}

Under regular perturbations as in \cite{LMSW96}, center-stable, 
center-unstable, and center manifolds, and Fenichel stable and 
unstable
fibers persist as in the standard theory. Under the singular 
perturbation, what are the objects that persist? We start with the 
linear
operator $L_\epsilon$.

\subsubsection{The Spectrum of $L_\epsilon$}

The spectrum of $L_\epsilon$ consists of only point spectrum. The 
eigenvalues of $L_\epsilon$ are:
\begin{equation}\label{leev} \mu^\pm_k=-\epsilon (\alpha+k^2)\pm 
k\sqrt{4\omega^2-k^2},\quad (k=1,2,\ldots ).\end{equation}
When $\omega \in \left( \frac{1}{2},1\right)$, only $\mu^\pm_1$ are 
real, and $\mu^\pm_k$ are complex for $k>1$. When $\omega \in
\left(1,\frac{3}{2}\right)$, only $\mu^\pm_1$ and $\mu^\pm_2$ are 
real, and $\mu^\pm_k$ are complex for $k>2$.

\setlength{\unitlength}{1in}
\begin{figure}
\begin{picture}(6.5,2.5)
\put(0,1.25){\line(2,0){2.5}}
\put(3.5,1.25){\line(2,0){2.5}}
\put(1.25,0){\line(0,1){2.5}}
\put(4.75,0){\line(0,1){2.5}}
\put(1.75,.5){$\epsilon =0$}
\put(5.25,.5){$\epsilon >0$}
\put(1.375,2.375){$\mu$}
\put(4.875,2.375){$\mu$}
\put(.5,1.21){$\bullet$}
\put(2,1.21){$\bullet$}
\put(4,1.21){$\bullet$}
\put(5.5,1.21){$\bullet$}
\put(1.21,.25){$\bullet$}
\put(1.21,.5){$\bullet$}
\put(1.21,.75){$\bullet$}
\put(1.21,1){$\bullet$}
\put(1.21,1.5){$\bullet$}
\put(1.21,1.75){$\bullet$}
\put(1.21,2){$\bullet$}
\put(1.21,2.25){$\bullet$}
\put(4,.25){$\bullet$}
\put(4.125,.5){$\bullet$}
\put(4.25,.75){$\bullet$}
\put(4.375,1){$\bullet$}
\put(4,2.25){$\bullet$}
\put(4.125,2){$\bullet$}
\put(4.25,1.75){$\bullet$}
\put(4.375,1.5){$\bullet$}
\end{picture}
\caption{The point spectra of the linear operator $L_\e$.}
\label{sgcb}
\end{figure}
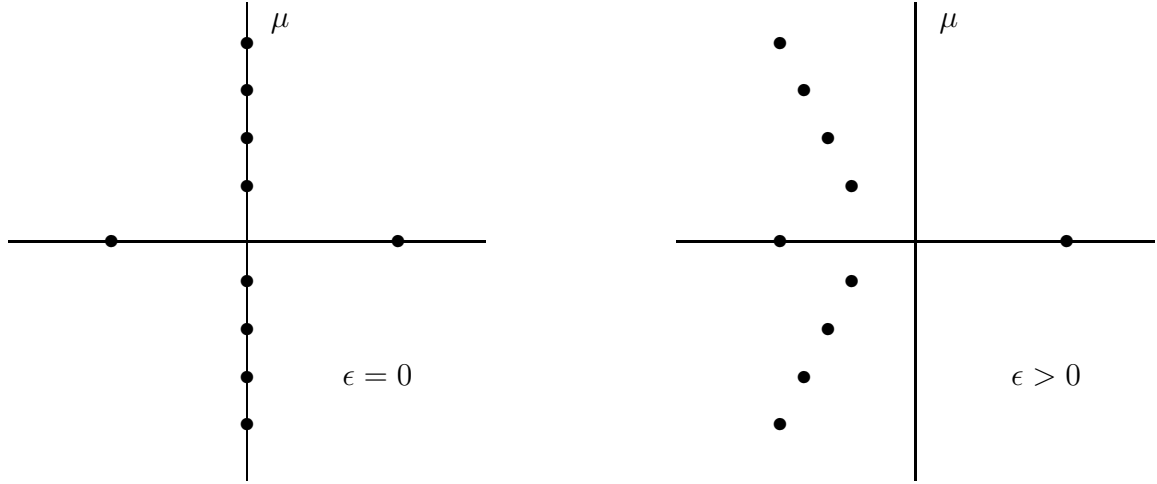

\begin{remark} The main difficulty introduced by the singular 
perturbation $\epsilon \partial^2_x f$ is the breaking of the 
spectral
gap condition. Figure~\ref{sgcb} shows the distributions of the 
eigenvalues when $\epsilon =0$ and $\epsilon \neq 0$. It clearly shows
the breaking of the stable spectral gap condition. As a result, 
center and center-unstable manifolds do not necessarily
persist. On the other hand, the unstable spectral gap condition is 
not broken. This gives the hope for the persistence of center-stable
manifold. Another case of persistence can be described as follows:
Notice that the plane $\Pi$ \eqref{Pi} is invariant 
under the PNLS flow \eqref{pnls}. When $\epsilon =0$, there is an
unstable fibration with base points in a neighborhood of the circle 
$S_\omega$ \eqref{rcl} in $\Pi$, as an invariant sub-fibration of the
unstable Fenichel fibration with base points in the center manifold. 
When $\epsilon >0$, the center manifold may not persist, but $\Pi$
persists, moreover, the unstable spectral gap condition is not 
broken, therefore, the unstable sub-fibration with base points in 
$\Pi$ may persist. This is the topics of this subsection. Since the 
semiflow generated by PNLS \eqref{pnls} is not a $C^1$ perturbation of
that generated by the unperturbed NLS due to the singular 
perturbation $\epsilon \partial^2_x$, standard results on 
persistence can not be applied.
\end{remark}
From now on, we will take the case of two unstable eigenvalues 
as our example to conduct the arguments. The case of one unstable 
eigenvalue is easier.
The eigenfunctions corresponding to the real 
eigenvalues are:
\begin{equation}\label{leef} e^\pm_k=e^{\pm i\vth_k}\cos kx,\quad 
e^{\pm i\vth_k}=\frac{k\mp
i\sqrt{4\omega^2-k^2}}{2\omega},\quad k=1,2.\end{equation}
Notice that they are independent of $\epsilon$.
The eigenspaces corresponding to the complex conjugate pairs of 
eigenvalues are given by:
\begin{equation*}E_k=\Span _{\mathbb{C}}\{ \cos k x \}.\end{equation*}
and have real dimension $2$.

\subsubsection{The Set-Up of Equations}

For the goal of this subsection, we need to single out the 
eigen-directions \eqref{leef}. Let
\begin{equation*}g=\sum_{\pm ;k=1,2}\xi_k^\pm e^\pm_k+h,\end{equation*}
where $\xi^\pm_k$ are real variables, and
\begin{equation*}\langle h\rangle=\langle h\cos x \rangle =\langle 
h\cos 2x \rangle =0.\end{equation*}
In terms of the coordinates $(\xi_k^\pm , J,\theta ,h)$, 
\eqref{nfe1}-\eqref{nfe3} can be rewritten as:
\begin{align}\label{eig1} 
\dot{\xi}^+_k&=\mu^+_k\xi^+_k+V^+_k\xi^+_k+\mathcal{N}^+_k,\quad 
(k=1,2),\\
\label{eig2} \dot{J}&=\epsilon \left[ -2\alpha (J+\omega^2)+2\beta 
\sqrt{J+\omega^2}\cos \theta \right] +\epsilon \mathcal{R}^J_2,\\
\label{eig3} \dot{\theta} & = -2J-\epsilon \beta \frac{\sin 
\theta}{\sqrt{J+\omega^2}}+\mathcal{R}^\theta_2,\\
\label{eig4} h_t&=L_\epsilon h+V_\epsilon h+\tilde{\mathcal{N}},\\
\label{eig5} 
\dot{\xi}^-_k&=\mu^-_k\xi^-_k+V^-_k\xi^-_k+\mathcal{N}^-_k,\quad 
(k=1,2),\end{align}
where $\mu^\pm_k$ are given in \eqref{leev}, $\mathcal{N}^\pm_k$ and 
$\tilde{\mathcal{N}}$ are projections of $\mathcal{N}$ to the
corresponding directions, and
\begin{align*}V^+_k\xi^+_k&=2c_kJ(\xi^+_k+\xi^-_k)+\epsilon \beta 
\frac{\sin \theta}{\sqrt{J+\omega^2}}(c^+_k\xi^+_k-c^-_k\xi^-_k),\\
V^-_k\xi^-_k&=-2c_kJ(\xi^+_k+\xi^-_k)+\epsilon \beta \frac{\sin 
\theta}{\sqrt{J+\omega^2}}(c_k^-\xi^+_k-c^+_k\xi^-_k),\\
c_k&=\frac{k}{\sqrt{4\omega^2-k^2}},\quad 
c^+_k=\frac{2\omega^2-k^2}{k\sqrt{4\omega^2-k^2}},\quad
c^-_k=\frac{2\omega^2}{k\sqrt{4\omega^2-k^2}}.\end{align*}

\subsubsection{Statement of the Unstable Fiber Theorem}

The main unstable fiber theorem can be stated as follows.
\begin{theorem}\label{UFT} There exists an annular neighborhood 
$\mathcal{A}$ of the circle $S_\omega$ \eqref{rcl} in $\Pi $
\eqref{Pi}, for any $p\in \mathcal{A}$, there is a local unstable 
fiber $\mathcal{F}^+_p$ which is a 2D surface. $\mathcal{F}^+_p$ has
the following properties:
\begin{enumerate}\item $\mathcal{F}^+_p$ is a $C^1$ smooth surface in 
$\|  \ \|_n$ norm, $\forall n\geq 1$.
\item $\mathcal{F}^+_p$ is also $C^1$ smooth in $\epsilon$, $\alpha$, 
$\beta$, $\omega$, and $p$ in $\|\ \|_n$ norm, for any $n\geq 1$,
$\epsilon \in [0,\epsilon_0)$ for some $\epsilon_0>0$.
\item $p\in \mathcal{F}^+_p$, $\mathcal{F}^+_p$ is tangent to $\Span 
\{ e^+_1,e^+_2\}$ at $p$ when $\epsilon=0$, where $e^+_k\ (k=1,2)$
are defined in \eqref{leef}.
\item $\mathcal{F}^+_p$ has the exponential decay property: Let $S^t$ 
be the evolution operator of \eqref{eig1}-\eqref{eig5}, $\forall
p_1\in \mathcal{F}^+_p$,
\begin{equation*}\frac{\| S^tp_1-S^tp\|_n}{\| p_1-p\|_n}\leq 
Ce^{\frac{\mu^+}{3}t},\quad \forall t\leq 0,\end{equation*}
where $\mu^+=\min \{ \mu^+_1,\mu^+_2\}$.
\item $\{ \mathcal{F}^+_p\}_{p\in \mathcal{A}}$ forms an invariant 
family of unstable fibers,
\begin{equation*}S^t\mathcal{F}^+_p\subset \mathcal{F}^+_{S^tp}\ ,\quad 
\forall t\in [-T,0],\end{equation*}
and $\forall T>0$ ($T$ can be $+\infty$), such that $S^tp\in 
\mathcal{A}$, $\forall t\in [-T,0]$.
\end{enumerate}\end{theorem}

\subsubsection{Proof of the Unstable Fiber Theorem}

There are two main approaches in establishing invariant manifolds and 
fibrations: 1. Lyapunov-Perron's method \cite{Per30}, 2. Hadamard's
method \cite{Had01}. Here we will adopt the Lyapunov-Perron's method, 
pay special attention to non-standard applications of the method, and
focus on the difficulties generated by the singular perturbation 
$\epsilon \partial^2_x$.

\begin{definition} For any $\delta >0$, we define the annular 
neighborhood of the circle $S_\omega$ \eqref{rcl} as
\begin{equation*}\mathcal{A}(\delta)=\{ (J,\theta )\mid \ |J|<\delta 
\}.\end{equation*}
\end{definition}
\begin{figure}
\vspace{1.5in}
\caption{The bump function $\eta$.}
\label{bump}
\end{figure}
\nid
To apply the Lyapunov-Perron's method, it is standard and necessary to 
modify the $J$ equation so that $\mathcal{A}(4\delta)$ is
overflowing invariant. Let $\eta \in C^\infty (R,R)$ be a ``bump" function:
\begin{equation*}\eta=\begin{cases} 0, & \text{in } (-2,2)\cup 
(-\infty,-6)\cup (6,\infty),\\
1, & \text{in } (3,5),\\
-1, & \text{in } (-5,-3),\end{cases}\end{equation*}
as shown in Figure \ref{bump}, $|\eta'|\leq 2$, $|\eta''|\leq C$. We 
modify the $J$ equation \eqref{eig2} as follows:
\begin{equation}\label{meig2}\dot{J} =\epsilon b\eta ( 
J/\delta)+\epsilon \left[ -2\alpha (J+\omega^2)+2\beta
\sqrt{J+\omega^2}\cos
\theta \right] +\epsilon \mathcal{R}^J_2,\end{equation}
where $b>2 (2\alpha \omega^2+2\beta \omega)$. Then 
$\mathcal{A}(4\delta)$ is overflowing invariant. There are two main 
points in
adopting the bump function:
\begin{enumerate}\item One needs $\mathcal{A}(4\delta)$ to be 
overflowing invariant so that a Lyapunov-Perron type integral
equation can be set up along orbits in $\mathcal{A}(4\delta)$ for 
$t\in (-\infty ,0)$.
\item One needs the vector field inside $\mathcal{A}(2\delta)$ to be 
unchanged so that results for the modified system can be claimed
for the original system in $\mathcal{A}(\delta)$.\end{enumerate}

\begin{remark} Due to the singular perturbation, the real part of 
$\mu^\pm_k$ approaches $-\infty$ as $k\to \infty$. Thus the $h$
equation \eqref{eig4} can not be modified to give overflowing flow. 
This rules out the construction of unstable fibers with base points
having general $h$ coordinates.\end{remark}
\nid
For any $(J_0,\theta_0)\in \mathcal{A}(4\delta)$, let
\begin{equation}\label{borb} J=J_*(t),\quad \theta=\theta_*(t), \quad 
t\in (-\infty ,0],\end{equation}
be the backward orbit of the modified system \eqref{meig2} and 
\eqref{eig3} with the initial point $(J_0,\theta_0)$. If
\begin{equation*}(\xi^+_k(t),J_*(t)+\tilde{J}(t),\theta_*(t)+\tilde{\theta}(t),h(t),\xi^-_k(t))\end{equation*}
is a solution of the modified full system, then one has
\begin{eqnarray}\dot{\xi}^+_k &=& \mu^+_k\xi^+_k+F^+_k,\quad (k=1,2)
\label{meq1}\\
u_t &=& Au+F,\label{meq2} 
\end{eqnarray}
where
\begin{align*} u&=\begin{pmatrix} \tilde{J} \\ \tilde{\theta}\\ h\\ 
\xi^-_1\\ \xi^-_2\end{pmatrix}, \quad A=\begin{pmatrix} 0 & 0& 0 & 0
& 0\\
-2 & 0 & 0 & 0 & 0\\ 0 & 0 & L_\epsilon & 0 & 0\\
0 & 0 & 0 & \mu^-_1 & 0\\
0 & 0 & 0 & 0 & \mu^-_2\end{pmatrix},\quad F=\begin{pmatrix} F_J\\ 
F_\theta\\ F_h\\
F^-_1\\ F^-_2\end{pmatrix},\\
&F^+_k  = V^+_k\xi^+_k+\mathcal{N}^+_k,\\
&F_J=\epsilon b\left[ \eta (J/\delta)-\eta 
(J_*(t)/\delta)\right]+\epsilon \left[
-2\alpha \tilde{J} +2\beta \sqrt{J+\omega^2}\cos \theta \right.\\
&\quad -\left. 2\beta \sqrt{J_*(t)+\omega^2}\cos \theta_*(t)\right] +\epsilon
\mathcal{R}^J_2,\\
&F_\theta = -\epsilon \beta \frac{\sin 
\theta}{\sqrt{J+\omega^2}}+\epsilon \beta \frac{\sin
\theta_*(t)}{\sqrt{J_*(t)+\omega^2}}+\mathcal{R}^\theta_2,\\
&F_h=V_\epsilon h+\tilde{\mathcal{N}},\\
&F_k^-=V^-_k\xi^-_k+\mathcal{N}^-_k,\\
&J=J_*(t)+\tilde{J},\quad \theta=\theta_*(t)+\tilde{\theta}.
\end{align*}
\nid
System \eqref{meq1}-\eqref{meq2} can be written in the equivalent 
integral equation form:
\begin{align}\label{eit1} \xi^+_k(t)&= 
\xi^+_k(t_0)e^{\mu^+_k(t-t_0)}+\int^t_{t_0}e^{\mu^+_k(t-\tau)}F^+_k(\tau 
)d\tau,\ \ (k=1,2)\\
\label{eit2} u(t)&= 
e^{A(t-t_0)}u(t_0)+\int^t_{t_0}e^{A(t-\tau)}F(\tau )d\tau.\end{align}
By virtue of the gap between $\mu^+_k$ and the real parts of the 
eigenvalues of $A$, one can introduce the following space: For $\sigma
\in \left( \frac{\mu^+}{100},\frac{\mu^+}{3}\right)$, $\mu^+=\min \{ 
\mu^+_1,\mu^+_2\}$, and $n\geq 1$, let
\begin{equation*}\begin{split} G_{\sigma ,n}&=\bigg \{ 
g(t)=(\xi^+_k(t),u(t))\bigg | \  t \in (-\infty ,0], \ g(t)\text{ is
continuous} \\ &\quad \text{in } t\text{ in } H^n\text{ norm }, 
\ \| g\|_{\sigma ,n}=\sup_{t\leq 0}e^{-\sigma t} [
\sum_{k=1,2}|\xi^+_k(t)|+\| u(t)\|_n ] <\infty\bigg \}\ . 
\end{split}\end{equation*}
$G_{\sigma ,n}$ is a Banach space under the norm $\| \cdot \|_{\sigma 
,n}$. Let $\mathcal{B}_{\sigma ,n}(r)$ denote the ball in
$G_{\sigma ,n}$ centered at the origin with radius $r$. Since $A$ 
only has point spectrum, the spectral mapping theorem is valid. It is
obvious that for $t\geq 0$,
\begin{equation*}\| e^{At}u\|_n\leq C(1+t)\| u\|_n,\end{equation*}
for some constant $C$. Thus, if $g(t)\in \mathcal{B}_{\sigma 
,n}(r)$, $r<\infty$ is a solution of \eqref{eit1}-\eqref{eit2}, by
letting $t_0\to -\infty$ in \eqref{eit2} and setting $t_0=0$ in 
\eqref{eit1}, one has
\begin{align}\label{per1} 
\xi^+_k(t)&=\xi^+_k(0)e^{\mu^+_kt}+\int^t_0e^{\mu^+_k(t-\tau)}F^+_k(\tau 
)d\tau,\ \ (k=1,2)\\
\label{per2} u(t)&=\int^t_{-\infty}e^{A(t-\tau)}F(\tau )d\tau.\end{align}
For $g(t)\in \mathcal{B}_{\sigma ,n}(r)$, let $\Gamma(g)$ be the 
map defined by the right hand side of \eqref{per1}-\eqref{per2}.
Then a solution of \eqref{per1}-\eqref{per2} is a fixed point of 
$\Gamma$. For any $n\geq 1$ and $\epsilon <\delta^2$, and $\delta$ and
$r$ are small enough, $F^+_k$ and $F$ are Lipschitz in $g$ with small 
Lipschitz constants. Standard arguments of the Lyapunov-Perron's
method readily imply the existence of a fixed point $g_*$ of $\Gamma$ 
in $\mathcal{B}_{\sigma ,n}(r)$. The difficulties lie in the
investigation on the regularity of $g_*$ with respect to $(\epsilon 
,\alpha ,\beta ,\omega ,J_0,\theta_0,\xi^+_k(0))$. That is our
focus. The most difficult one is the regularity with respect to 
$\epsilon$ due to the singular perturbation, which is our further 
focus.
Formally differentiating $g_*$ in \eqref{per1}-\eqref{per2} with 
respect to $\epsilon$, one gets
\begin{align}\begin{split}\label{dper1} \xi^+_{k,\epsilon}(t)&= 
\int^t_0e^{\mu^+_k(t-\tau)}\left[ \partial_uF^+_k\cdot u_\epsilon
+\sum_{\ell=1,2}\partial_{\xi_\ell^+}F^+_k\cdot 
\xi^+_{\ell,\epsilon}\right](\tau )d\tau\\
&\quad +\mathcal{R}^+_k(t),\quad (k=1,2)\end{split}\\
u_\epsilon (t)&=\int^t_{-\infty}e^{A(t-\tau)}\left[ \partial_uF\cdot 
u_\epsilon +\sum_{\ell=1,2}\partial_{\xi_\ell^+}F\cdot
\xi^+_{\ell,\epsilon}\right](\tau )d\tau 
+\mathcal{R}(t),\label{dper2}\end{align}
where
\begin{align}\begin{split}\label{wdp1}
\mathcal{R}^+_k(t)&=\xi^+_k(0)\mu^+_{k,\epsilon}te^{\mu^+_kt}+\int^t_0\mu^+_{k,\epsilon}(t-\tau)e^{\mu^+_k(t-\tau)}F^+_k(\tau 
)d\tau\\
&\quad + \int^t_0e^{\mu^+_k(t-\tau)}[\partial_\epsilon 
F^+_k+\partial_{u_*}F^+_k\cdot u_{*,\epsilon}](\tau )d\tau,\end{split}\\
\begin{split}\label{wdp2} \mathcal{R}(t)&= 
\int^t_{-\infty}(t-\tau)A_\epsilon e^{A(t-\tau)}F(\tau )d\tau\\
&\quad +\int^+_{-\infty}e^{A(t-\tau)}[\partial_\epsilon 
F+\partial_{u_*}F\cdot u_{*,\epsilon}](\tau)d\tau,\end{split}\\
\label{wdp3} \mu^+_{k,\epsilon}&=-(\alpha +k^2),\quad k=1,2,\\
A_\epsilon &=\begin{pmatrix} 0 & 0 & 0 & 0 & 0\\
0 & 0 & 0 & 0 & 0\\
0 & 0 & -\alpha+\partial^2_x & 0 & 0\\
0 & 0 & 0 & -(\alpha+1) & 0\\
0 & 0 & 0 & 0 & -(\alpha+4)\end{pmatrix},\\
u_*&=(J_*,\theta_*,0,0,0)^T,\label{wdp4}\end{align}
where $T=$\ transpose, and $(J_*,\theta_*)$ are given in \eqref{borb}. 
The troublesome terms are the ones containing $A_\epsilon $ or
$u_{*,\epsilon}$ in \eqref{wdp1}-\eqref{wdp2}.
\begin{equation} 
\| A_\epsilon F\|_n\leq C\| F\|_{n+2}\leq \tilde{c} 
( \|u\|_{n+2}+\sum_{k=1,2}|\xi^+_k|),
\label{upre}
\end{equation}
where $\tilde{c}$ is small when $(\ \cdot \ )$ on the right hand side 
is small.
\begin{align*}\begin{split} \partial_{J_*}F_J\cdot J_{*,\epsilon}&= 
\frac{\epsilon}{\delta}b[ \eta'(
J/\delta)-\eta'( J_*/\delta)] 
J_{*,\epsilon}\\
&\quad + \epsilon [ \beta \frac{\cos 
\theta}{\sqrt{J+\omega^2}}-\beta\frac{\cos 
\theta_*}{\sqrt{J_*+\omega^2}}]
J_{*,\epsilon}\\
&\quad +\text{easier terms}.\end{split}\\
\begin{split} |\partial_{J_*}F_J\cdot J_{*,\epsilon}|&\leq 
\frac{\epsilon}{\delta^2}b\sup_{0\leq \hat{\gamma}\leq 1}\left| \eta''(
[\hat{\gamma} J_*+(1-\hat{\gamma})J]/\delta)\right | \ |\tilde{J}| 
\ |J_{*,\epsilon}|\\
&\quad + \epsilon \beta 
C(|\tilde{J}|+|\tilde{\theta}|)|J_{*,\epsilon}|+\text{ easier terms}\\
& \leq C_1(|\tilde{J}|+|\tilde{\theta}|)|J_{*,\epsilon}|+\text{ 
easier terms,}\end{split}\\
\begin{split} \sup_{t\leq 0} e^{-\sigma t} |\partial_{J_*}F_J\cdot 
J_{*,\epsilon}|&\leq C_1\sup_{t\leq 0}[ e^{-(\sigma
+\tilde{\nu} )t}(|\tilde{J}|+|\tilde{\theta}|)]\sup_{t\leq 0}[e^{\tilde{\nu} 
t}|J_{*,\epsilon}|]\\
&\quad + \text{easier terms},\end{split}\end{align*}
where $\sup_{t\leq 0}e^{\tilde{\nu} t}|J_{*,\epsilon}|$ can be bounded when 
$\epsilon$ is sufficiently small for any fixed $\tilde{\nu} >0$, through a
routine estimate on Equations~\eqref{meig2} and \eqref{eig3} for 
$(J_*(t),\theta_*(t))$. Other terms involving $u_{*,\epsilon}$ can be
estimated similarly. Thus, the $\| \ \|_{\sigma ,n}$ norm of terms 
involving $u_{*,\epsilon}$ has to be bounded by $\|\ \|_{\sigma +\tilde{\nu}
,n}$ norms. This leads to the standard rate condition for the 
regularity of invariant manifolds. That is, the regularity is 
controlled by
the spectral gap. The $\| \ \|_{\sigma ,n}$ norm of the term involving 
$A_\epsilon$ has to be bounded by $\| \ \|_{\sigma ,n+2}$ norms. This
is a new phenomenon caused by the singular perturbation. This problem 
is resolved by virtue of a special property of the fixed point $g_*$
of $\Gamma$.
Notice that if $\sigma_2\geq \sigma_1$, $n_2\geq n_1$, then 
$G_{\sigma_2,n_2}\subset G_{\sigma_1,n_1}$. Thus by the uniqueness of 
the
fixed point, if $g_*$ is the fixed point of $\Gamma$ in 
$G_{\sigma_2,n_2}$, $g_*$ is also the fixed point of $\Gamma$ in
$G_{\sigma_1,n_1}$. Since $g_*$ exists in $G_{\sigma ,n}$ for an 
fixed $n\geq 1$ and $\sigma \in (
\frac{\mu^+}{100},\frac{\mu^+}{3}-10\tilde{\nu})$ where 
$\tilde{\nu}$ is small enough,
\begin{align*} \| \mathcal{R}^+_k\| _{\sigma,n}&\leq C_1+C_2\| 
g_*\|_{\sigma+\tilde{\nu},n},\\
\| \mathcal{R}\|_{\sigma ,n}&\leq C_3\| g_*\|_{\sigma ,n+2}+C_4\| 
g_*\|_{\sigma+\tilde{\nu},n}+C_5,\end{align*}
where $C_j\ (1\leq j\leq 5)$ depend upon $\| g_*(0)\|_n$ and $\| 
g_*(0)\|_{n+2}$. Let
\begin{equation*}M=2(\| \mathcal{R}^+_k\|_{\sigma 
,n}+\|\mathcal{R}\|_{\sigma, n}),\end{equation*}
and $\Gamma'$ denote the linear map defined by the right hand sides 
of \eqref{dper1} and \eqref{dper2}. Since the terms
$\partial_uF^+_k$, $\partial_{\xi^+_\ell}F^+_k$, $\partial_uF$, and 
$\partial_{\xi^+_\ell}F$ all have small $\|\ \|_n$ norms, $\Gamma '$
is a contraction map on $\mathcal{B}(M)\subset L(\mathcal{R}, G_{\sigma 
,n})$, where $\mathcal{B}(M)$ is the ball of radius $M$. Thus
$\Gamma '$ has a unique fixed point $g_{*,\epsilon}$. Next one needs 
to show that $g_{*,\epsilon}$ is indeed the partial derivative of
$g_*$ with respect to $\epsilon$. That is, one needs to show
\begin{equation}\label{ddef} \lim_{\Delta \epsilon \to 0}\frac{\| 
g_*(\epsilon +\Delta \epsilon)-g_*(\epsilon)-g_{*,\epsilon} \Delta
\epsilon \|_{\sigma ,n}}{\Delta \epsilon}=0.\end{equation}
\nid
This has to be accomplished directly from Equations 
\eqref{per1}-\eqref{per2}, \eqref{dper1}-\eqref{dper2} satisfied by 
$g_*$ and
$g_{*,\epsilon}$. The most troublesome estimate is still the one 
involving $A_\epsilon$. First, notice the fact that $e^{\epsilon
\partial^2_x}$ is holomorphic in $\epsilon$ when $\epsilon>0$, and 
not differentiable at $\epsilon=0$. Then, notice that $g_*\in
G_{\sigma ,n}$ for any $n\geq 1$, thus, $e^{\epsilon 
\partial_x^2}g_*$ is differentiable, up to certain order $m$, in 
$\epsilon$ at
$\epsilon=0$ from the right, i.e.
\begin{equation*}(d^+/d\epsilon)^me^{\epsilon 
\partial^2_x}g_*|_{\epsilon =0}\end{equation*}
exists in $H^n$. Let
\begin{equation*}\begin{split} z(t,\Delta \epsilon)&= e^{(\epsilon 
+\Delta \epsilon)t\partial^2_x}g_*-e^{\epsilon
t\partial^2_x}g_*-(\Delta \epsilon) t\partial ^2_x e^{\epsilon 
t\partial _x^2}g_*\\
&= e^{\epsilon t\partial^2_x}w(\Delta \epsilon),\end{split}\end{equation*}
where $t\geq 0$, $\Delta \epsilon >0$, and
\begin{equation*}w(\Delta \epsilon) =e^{(\Delta \epsilon ) t\partial 
_x^2}g_*-g_*-(\Delta \epsilon ) t\partial^2_x g_*.\end{equation*}
Since $w(0)=0$, by the Mean Value Theorem, one has
\begin{equation*}\| w(\Delta \epsilon )\|_n=\| w(\Delta \epsilon 
)-w(0)\|_n\leq \sup_{0\leq \lambda \leq 1}\| \frac{dw}{d\Delta
\epsilon}(\lambda
\Delta \epsilon)\|_n|\Delta \epsilon|,\end{equation*}
where at $\lambda=0$, $\frac{d}{d\Delta \epsilon}=\frac{d^+}{d\Delta 
\epsilon}$, and
\begin{equation*}\frac{dw}{d\Delta \epsilon}=t[e^{(\Delta \epsilon )
t\partial^2_x}\partial^2_x g_*-\partial^2_x g_*].
\end{equation*}
Since $\frac{dw}{d\Delta \epsilon}(0)=0$, by the Mean Value Theorem 
again, one has
\[
\| \frac {dw}{d\Dl \e }(\la \Dl \e )\|_n = 
\| \frac {dw}{d\Dl \e }(\la \Dl \e ) - 
\frac {dw}{d\Dl \e }(0)\|_n \leq \sup_{0 \leq \la_1 \leq 1}
\| \frac {d^2w}{d\Dl \e^2 }(\la_1 \la \Dl \e )\|_n |\la | |\Dl \e |\ ,
\]
where
\[
\frac {d^2w}{d\Dl \e^2 } = t^2 [e^{(\Dl \e )t \pa_x^2} \pa_x^4 g_* ]\ .
\]
Therefore, one has the estimate
\begin{equation}
\| z(t,\Dl \e )\|_n \leq |\Dl \e |^2 t^2 \| g_* \|_{n+4} \ .
\label{defes}
\end{equation}
This estimate is sufficient for handling the estimate involving 
$A_\e$. The estimate involving $u_{*,\e}$ can be handled in a similar 
manner. For instance, let 
\[
\tilde{z}(t,\Dl \e ) = F(u_*(t,\e +\Dl \e ))- F(u_*(t,\e ))
-\Dl \e \pa_{u_*}F \cdot u_{*,\e}\ ,
\]
then
\[
\| \tilde{z}(t,\Dl \e )\|_{\sg , n} \leq |\Dl \e |^2 
\sup_{0 \leq \la \leq 1} \| [u_{*,\e} \cdot \pa^2_{u_*}F \cdot u_{*,\e}
+\pa_{u_*}F \cdot u_{*,\e \e}](\la \Dl \e )\|_{\sg , n} \ .
\]
From the expression of $F$ (\ref{meq2}), one has
\begin{equation*}\begin{split} & \| u_{*,\epsilon}\cdot 
\partial^2_{u_*}F\cdot u_{*,\epsilon}+\partial_{u_*}F\cdot
u_{*,\epsilon\epsilon}\|_{\sigma,n}\\
&\quad \leq C_1\| g_*\|_{\sigma +2\tilde{\nu} ,n} [( \sup_{t\leq 
0}e^{\tilde{\nu} t}|u_{*,\epsilon}|)^2+\sup_{t\leq 0}e^{2\tilde{\nu}
t}|u_{*,\epsilon \epsilon}|],\end{split}\end{equation*}
and the term $[ \  ]$ on the right hand side can be easily shown to be 
bounded. In conclusion, let
\begin{equation*}h=g_*(\epsilon +\Delta 
\epsilon)-g_*(\epsilon)-g_{*,\epsilon}\Delta 
\epsilon,\end{equation*}
one has the estimate
\begin{equation*}\| h\|_{\sigma ,n}\leq \tilde{\k} \| 
h\|_{\sigma,n}+|\Delta \epsilon|^2\tilde{C}(\| g_*\|_{\sigma ,n+4};\|
g_*\|_{\sigma +2\tilde{\nu} ,n}),\end{equation*}
where $\tilde{\k}$ is small, thus
\begin{equation*}\| h\|_{\sigma ,n}\leq 2 |\Delta \epsilon 
|^2\tilde{C}(\| g_*\|_{\sigma,n+4};\| g_*\|_{\sigma +2\tilde{\nu}
,n}).\end{equation*} This implies that
\begin{equation*}\lim_{\Delta \epsilon \to 0}\frac{\| h\|_{\sigma 
,n}}{|\Delta \epsilon|}=0,\end{equation*}
which is \eqref{ddef}.

Let $g_*(t)=(\xi^+_k(t),u(t))$. First, let me comment on $\left. 
\frac{\partial u}{\partial
\xi^+_k(0)}\right|_{\xi^+_k(0)=0,\epsilon =0}=0$. From \eqref{per2}, one has
\begin{equation*}\| \frac{\partial u}{\partial 
\xi^+_k(0)}\| _{\sigma, n}\leq \k_1\| \frac{\partial 
g_*}{\partial \xi^+_k(0)}\|_{\sigma ,n},\end{equation*}
by letting $\xi^+_k(0)\to 0$ and $\epsilon \to 0^+$, $\k_1\to 0$. Thus
\begin{equation*}\left. \frac{\partial u}{\partial 
\xi^+_k(0)}\right|_{\xi^+_k(0)=0,\epsilon =0}=0.\end{equation*}
I shall also comment on ``exponential decay" property. Since $\| 
g_*\|_{\frac{\mu^+}{3},n}\leq r$,
\begin{equation*}\| g_*(t)\|_n\leq re^{\frac{\mu^+}{3}t},\quad 
\forall t\leq 0.\end{equation*}
\begin{definition} Let $g_*(t)=(\xi^+_k(t), u(t))$, where
\begin{equation*}u(0)=\int^0_{-\infty}e^{A(t-\tau)}F(\tau)d\tau\end{equation*}
depends upon $\xi^+_k(0)$. Thus
\begin{equation*}u^0_*:\xi^+_k(0)\mapsto u(0),\end{equation*}
defines a $2D$ surface, which we call an unstable fiber denoted by 
$\mathcal{F}^+_p$, where $p=(J_0,\theta_0)$ is the base point,
$\xi^+_k(0)\in [-r,r]\times [-r,r]$.
\end{definition}
Let $S^t$ denote the evolution operator of \eqref{meq1}-\eqref{meq2}, then
\begin{equation*}S^t\mathcal{F}_p^t\subset \mathcal{F}^t_{S^tp},\quad 
\forall t\leq 0.\end{equation*}
That is, $\{ \mathcal{F}^+_p\}_{p\in \mathcal{A}(4\delta)}$ is an 
invariant family of unstable fibers. The proof of the Unstable Fiber
Theorem is finished.\hfill$\Box$

\begin{remark} If one replaces the base orbit $(J_*(t),\theta_*(t))$ 
by a general orbit for which only $\|\ \|_n$ norm is bounded, then
the estimate \eqref{upre} will not be possible. The $\|\ \|_{\sigma 
,n+2}$ norm of the fixed point $g_*$ will not be bounded either. In
such case, $g_*$ may not be smooth in $\epsilon$ due to the singular 
perturbation.\end{remark}

\begin{remark} Smoothness of $g_*$ in $\epsilon$ at $\epsilon =0$ is 
the key point of the entire argument in this article. In the global
theory in later sections, information is known at $\epsilon =0$. This 
key point will 
link ``$\epsilon =0$" information to ``$\epsilon \neq 0$" studies.
Only continuity in $\epsilon $ at $\epsilon =0$ is not enough for the 
study. The beauty of the entire theory is reflected by the fact
that although $e^{\epsilon \partial^2_x}$ is not holomorphic at 
$\epsilon =0$, $e^{\epsilon \partial^2_x}g_*$ can be smooth at
$\epsilon =0$ from the right, up to certain order depending upon the 
regularity of $g_*$. This is the beauty of the singular
perturbation.
\end{remark}

\subsection{Center-Stable Manifold}

We start with Equations \eqref{eig1}-\eqref{eig5}, let
\begin{equation}\label{2cr} v=\begin{pmatrix} J\\ \theta \\ h \\ 
\xi^-_1\\ \xi^-_2\end{pmatrix},\quad \tilde{v}=\begin{pmatrix}
J\\ h\\ \xi^-_1\\ \xi^-_2\end{pmatrix},\end{equation}
and let $E_n(r)$ be the tubular neighborhood of $S_\om$ \eqref{rcl}:
\begin{equation}\label{defEn} E_n(r)=\{ (J,\theta 
,h,\xi^-_1,\xi^-_2)\in H^n\mid \ \| \tilde{v}\|_n\leq r\}.\end{equation}
$E_n(r)$ is of codimension $2$ in the entire phase space 
coordinatized by $(\xi^+_1,\xi^+_2,J,\theta
,h,\xi^-_1,\xi^-_2)$.

\subsubsection{Statement of the Center-Stable Manifold Theorem}

\begin{theorem}\label{CSM} There exists a $C^1$ smooth codimension 2 
locally invariant center-stable manifold $W^{cs}_n$ in $H^n$ for any
$n\geq 1$.
$W^{cs}_n$ can be represented as the graph of a $C^1$ function 
$\xi^+_*:E_n(r)\to R^2$, for some $r>0$.
\begin{enumerate}\item At points in the subset $W^{cs}_{n+4}$ of 
$W^{cs}_n$, $W^{cs}_n$ is $C^1$ smooth in $\epsilon$ for $\epsilon \in
[0,\epsilon_0)$ and some $\epsilon_0 >0$. That is, if $v\in 
E_{n+4}(r)\subset E_n(r)$, then $\xi^+_*(v)$ is $C^1$ smooth in 
$\epsilon$, in
$H^n$ norm. Moreover, $\partial_\epsilon \xi^+_*(v)$ is uniformly 
bounded in $v\in E_{n+4}(r)$ and $\epsilon \in [0,\epsilon_0)$.
\item $W^{cs}_n$ is $C^1$ smooth in $(\alpha ,\beta ,\omega)$.
\item The annular neighborhood $\mathcal{A}$ in Theorem~\ref{UFT} is 
included in $W^{cs}_n$, i.e.
$\xi^+_*(J,\theta, 0,0,0)=0$. Along the circle $S_\omega$ 
\eqref{rcl}, $W^{cs}_n$ is tangent to $E_n(r)$ when
$\epsilon=0$, i.e. $\partial_v\xi^+_*(0,\theta,0,0)=0$ when 
$\epsilon=0$. $W^{cs}_n$ is $C^1$ close to
$E_n(r)$, i.e. $\| \partial_v\xi^+_*\| \leq Cr$.\end{enumerate}\end{theorem}

\begin{remark}\label{csnr} $C^1$ regularity in $\epsilon$ is 
crucial in locating a homoclinic orbit. As can be seen later, one
has detailed information on certain unperturbed (i.e. $\epsilon=0$) 
homoclinic orbit, which will be used in tracking candidates for a
perturbed homoclinic orbit. In particular, Melnikov measurement will 
be needed. Melnikov measurement measures zeros of
$\mathcal{O}(\epsilon)$ signed distances, thus, the perturbed orbit 
needs to be $\mathcal{O}(\epsilon)$ close to the unperturbed orbit in
order to perform Melnikov measurement.\end{remark}

\subsubsection{Proof of the Center-Stable Manifold Theorem}

Let $\chi\in C^\infty(R,R)$ be a ``cut-off" function:
\begin{equation*}\chi=\begin{cases} 0, & \text{in } (-\infty ,-4)\cup 
(4,\infty),\\
1, & \text{in } (-2,2).\end{cases}\end{equation*}
We apply the cut-off \begin{equation*}\chi_\delta =\chi ( 
\|\tilde{v}\|_n/\delta) \chi (
\xi^+_1/\delta)\chi (
\xi^+_2/\delta)\end{equation*} to 
Equations~\eqref{eig1}-\eqref{eig5}, so that the equations in a 
tubular neighborhood
of the circle
$S_\omega$ \eqref{rcl} are unchanged, and linear outside a 
bigger tubular neighborhood. The modified equations take the form:
\begin{align}\label{ceq1} 
\dot{\xi}_k^+&=\mu^+_k\xi^+_k+\tilde{F}^+_k, \quad (k=1,2)\\
\label{ceq2}v_t&=Av+\tilde{F},\end{align}
where $A$ is given in \eqref{meq2},
\begin{align*} \tilde{F}^+_k&=\chi_\delta[V^+_k\xi^+_k+\mathcal{N}^+_k],\\
\tilde{F}&=(\tilde{F}_J,\tilde{F}_\theta 
,\tilde{F}_h,\tilde{F}^-_1,\tilde{F}^-_2)^T,\quad 
T=\text{transpose},\\
\tilde{F}_J&=\chi_\delta \ \e \left[-2\alpha (J+\omega^2)+2\beta 
\sqrt{J+\omega^2}\cos \theta+\mathcal{R}^J_2\right],\\
\tilde{F}_\theta&=\chi _\delta \left[ -\epsilon \beta \frac{\sin 
\theta}{\sqrt{J+\omega^2}}+\mathcal{R}^\theta_2\right],\\
\tilde{F}_h&=\chi_\delta [V_\epsilon h+\tilde{\mathcal{N}}],
\ \ (k=1,2)\\
\tilde{F}_k^-&=\chi_\delta [V^-_k\xi^-_k+ \mathcal{N}^-_k],\end{align*}
Equations~\eqref{ceq1}-\eqref{ceq2} can be written in the equivalent 
integral equation form:
\begin{align}\label{csit1} 
\xi^+_k(t)&=\xi^+_k(t_0)e^{\mu^+_k(t-t_0)}+\int^t_{t_0}e^{\mu^+_k(t-\tau 
)}\tilde{F}^+_k(\tau 
)d\tau,\\
\label{csit2} 
v(t)&=e^{A(t-t_0)}v(t_0)+\int^t_{t_0}e^{A(t-\tau)}\tilde{F}(\tau 
)d\tau.\end{align}
We introduce the following space: For $\sigma \in \left( 
\frac{\mu^+}{100},\frac{\mu^+}{3}\right)$, $\mu^+=\min \{ 
\mu^+_1,\mu^+_2\}$,
and $n\geq 1$, let
\begin{equation*}\begin{split} \tilde{G}_{\sigma ,n}&= \bigg \{ 
g(t)=(\xi^+_k(t),v(t)) \bigg | \  t\in [0,\infty),g(t)\text{ is continuous 
in } t\\
&\quad \text{in } H^n\norm ,\| g\|_{\sigma ,n}=\sup_{t\geq 
0}e^{-\sigma t}[ \sum_{k=1,2}|\xi^+_k(t)|+\| v(t)\|_n ]
<\infty \bigg \}\ .\end{split}\end{equation*}
$\tilde{G}_{\sigma,n}$ is a Banach space under the norm $\| \cdot 
\|_{\sigma ,n}$. Let $\tilde{\mathcal{A}}_{\sigma ,n}(r)$ denote the
closed tubular neighborhood of $S_\omega$ \eqref{rcl}:
\begin{equation*} \tilde{\mathcal{A}}_{\sigma ,n}(r)= \bigg \{ g(t) 
=(\xi^+_k(t),v(t))\in \tilde{G}_{\sigma,n}\bigg | \ \sup_{t\geq
0}e^{-\sigma t}[ \sum_{k=1,2}|\xi^+_k(t)|+\| 
\tilde{v}(t)\|_n] \leq r\bigg \}\ ,\end{equation*}
where $\tilde{v}$ is defined in \eqref{2cr}. If $g(t)\in 
\tilde{\mathcal{A}}_{\sigma ,n}(r)$, $r<\infty$, is a solution of
\eqref{csit1}-\eqref{csit2}, by letting $t_0\to +\infty$ in 
\eqref{csit1} and setting $t_0=0$ in \eqref{csit2}, one has
\begin{align}\label{cspe1} 
\xi^+_k(t)&=\int^t_{+\infty}e^{\mu^+_k(t-\tau)}\tilde{F}^+_k(\tau 
)d\tau, \quad (k=1,2)\\
\label{cspe2}v(t)&=e^{At}v(0)+\int^t_0e^{A(t-\tau)}\tilde{F}(\tau 
)d\tau .\end{align}
For any $g(t)\in \tilde{\mathcal{A}}_{\sigma ,n}(r)$, let 
$\tilde{\Gamma}(g)$ be the map defined by the right hand side of
\eqref{cspe1}-\eqref{cspe2}. In contrast to the map $\Gamma$ defined 
in \eqref{per1}-\eqref{per2}, $\tilde{\Gamma}$ contains constant
terms of order $\mathcal{O}(\epsilon)$, e.g. $\tilde{F}_J$ and 
$\tilde{F}_\theta$ both contain such terms. Also,
$\tilde{\mathcal{A}}_{\sigma, n}(r)$ is a tubular neighborhood of the 
circle $S_\omega$ \eqref{rcl} instead of the ball
$\mathcal{B}_{\sigma ,n}(r)$ for $\Gamma$. Fortunately, these facts 
will not create any difficulty in showing $\tilde{\Gamma}$ is a
contraction on $\tilde{\mathcal{A}}_{\sigma ,n}(r)$. For any $n\geq 
1$ and $\epsilon <\delta^2$, and $\delta $ and $r$ are small enough,
$\tilde{F}_k^+$ and $\tilde{F}$ are Lipschitz in $g$ with small 
Lipschitz constants. $\tilde{\Gamma}$ has a unique fixed point
$\tilde{g}_*$ in $\tilde{\mathcal{A}}_{\sigma ,n}(r)$, following from 
standard arguments. For the regularity of $\tilde{g}_*$ with
respect to $(\epsilon, \alpha ,\beta ,\omega ,v(0))$, the most 
difficult one is of course with respect to $\epsilon$ due to the 
singular
perturbation. Formally differentiating $\tilde{g}_*$ in 
\eqref{cspe1}-\eqref{cspe2} with respect to $\epsilon$, one gets
\begin{align}\begin{split}\label{dcsp1} 
\xi^+_{k,\epsilon}(t)&=\int^t_{+\infty}e^{\mu^+_k(t-\tau)}\left[
\sum_{\ell=1,2}\partial_{\xi^+_\ell}\tilde{F}^+_k\cdot 
\xi^+_{\ell,\epsilon}+\partial_v\tilde{F}^+_k\cdot v_\epsilon 
\right](\tau
)d\tau\\
  & \quad +\tilde{R}^+_k(t),\quad (k=1,2)\end{split}\\
\label{dcsp2}v(t)&=\int^t_0e^{A(t-\tau)}\left[ 
\sum_{\ell=1,2}\partial_{\xi^+_\ell}\tilde{F}\cdot
\xi^+_{\ell,\epsilon}+\partial_v\tilde{F}\cdot v_\epsilon 
\right](\tau )d\tau +\tilde{R}(t),\end{align}
where
\begin{align}\begin{split}\label{wdcs1} 
\tilde{R}_k^+(t)&=\int^t_{+\infty}\mu^+_{k,\epsilon}(t-\tau)
e^{\mu^+_k(t-\tau)}\tilde{F}_k^+(\tau )d\tau
+\int^t_{+\infty}e^{\mu^+_k(t-\tau)}\partial _\epsilon 
\tilde{F}^+_k(\tau )d\tau,\end{split}\\
\begin{split}\label{wdcs2} \tilde{R}(t)&= tA_\epsilon 
e^{At}v(0)+\int^t_0(t-\tau)A_\epsilon e^{A(t-\tau)}\tilde{F}(\tau 
)d\tau + \int^t_0e^{A(t-\tau )}\partial_\epsilon \tilde{F}(\tau 
)d\tau,\end{split}\end{align}
and $\mu^+_{k,\epsilon}$ and $A_\epsilon$ are given in 
\eqref{wdp3}-\eqref{wdp4}. The troublesome terms are the ones 
containing
$A_\epsilon$ in \eqref{wdcs2}. These terms can be handled in the same 
way as in the Proof of the Unstable Fiber Theorem. The crucial
fact utilized is that if $v(0)\in H^{n_1}$, then $\tilde{g}_*$ is the 
unique fixed point of $\tilde{\Gamma}$ in both
$\tilde{G}_{\sigma,n_1}$ and $\tilde{G}_{\sigma ,n_2}$ for any $n_2\leq n_1$.

\begin{remark} In the Proof of the Unstable Fiber Theorem, the 
arbitrary initial data in \eqref{per1}-\eqref{per2} are $\xi^+_k(0)\ $
$(k=1,2)$ which are scalars. Here the arbitrary initial datum in 
\eqref{cspe1}-\eqref{cspe2} is $v(0)$ which is a function of $x$. 
If $v(0)\in H^{n_2}$ but not $H^{n_1}$ for some $n_1>n_2$, 
then $\tilde{g}_*\notin \tilde{G}_{\sigma ,n_1}$, in contrast to the 
case of
\eqref{per1}-\eqref{per2} where $g_*\in G_{\sigma ,n}$ for any fixed 
$n\geq 1$. The center-stable manifold $W^{cs}_n$ stated in the
Center-Stable Manifold Theorem will be defined through $v(0)$. This 
already illustrates why $W^{cs}_n$ has the regularity in $\epsilon$
as stated in the theorem.\end{remark}
\nid
We have
\begin{align*} \| \tilde{R}^+_k\|_{\sigma ,n}&\leq \tilde{C_1},\\
\| \tilde{R}\|_{\sigma ,n}&\leq \tilde{C}_2\| \tilde{g}_*\|_{\sigma 
,n+2}+\tilde{C}_3,\end{align*}
for $\tilde{g}_*\in \tilde{\mathcal{A}}_{\sigma , n+2}(r)$, where 
$\tilde{C}_j\ (j=1,2,3)$ are constants depending in particular upon the
cut-off in $\tilde{F}^+_k$ and $\tilde{F}$. Let $\tilde{\Gamma}'$ 
denote the linear map defined by the right hand sides of
\eqref{dcsp1}-\eqref{dcsp2}. If $v(0)\in H^{n+2}$ and $\tilde{g}_*\in 
\tilde{\mathcal{A}}_{\sigma,n+2}(r)$, standard argument shows that
$\tilde{\Gamma}'$ is a contraction map  on a closed ball in 
$L(R,\tilde{G}_{\sigma,n})$. Thus $\tilde{\Gamma}'$ has a unique 
fixed point
$\tilde{g}_{*,\epsilon}$. Furthermore, if $v(0)\in H^{n+4}$ and 
$\tilde{g}_*\in \tilde{\mathcal{A}}_{\sigma ,n+4}(r)$, one has that
$\tilde{g}_{*,\epsilon}$ is indeed the derivative of $\tilde{g}_*$ in 
$\epsilon$, following the same argument as in the Proof of the
Unstable Fiber Theorem. Here one may be able to replace the 
requirement $v(0)\in H^{n+4}$ and $\tilde{g}_*\in
\tilde{\mathcal{A}}_{\sigma ,n+4}(r)$ by just $v(0)\in H^{n+2}$ and 
$\tilde{g}_*\in \tilde{\mathcal{A}}_{\sigma ,n+2}(r)$. But we are
not interested in sharper results, and the current result is 
sufficient for our purpose.

\begin{definition} For any $v(0)\in E_n(r)$ where $r$ is sufficiently 
small and $E_n(r)$ is defined in \eqref{defEn}, let
$\tilde{g}_*(t)=(\xi^+_k(t),v(t))$ be the fixed point of 
$\tilde{\Gamma}$ in $\tilde{G}_{\sigma ,n}$, where one has
\begin{equation*}\xi^+_k(0)=\int^0_{+\infty}e^{\mu^+_k(t-\tau)}\tilde{F}^+_k(\tau 
)d\tau ,\quad (k=1,2)\end{equation*}
which depend upon $v(0)$. Thus
\begin{equation*} \xi^+_*:v(0)\mapsto \xi^+_k(0),\quad (k=1,2)\end{equation*}
defines a codimension $2$ surface, which we call center-stable 
manifold denoted by $W^{cs}_n$.\end{definition}

The regularity of the fixed point $\tilde{g}_*$ immediately implies 
the regularity of $W^{cs}_n$. We have sketched the proof of the most
difficult regularity, i.e. with respect to $\epsilon$. Uniform 
boundedness of $\partial_\epsilon \xi^+_*$ in $v(0)\in E_{n+4}(r)$ and
$\epsilon \in [0,\epsilon _0)$, is obvious. Other parts of the 
detailed proof is completely standard. We have that $W^{cs}_n$ is a 
$C^1$
locally invariant submanifold which is $C^1$ in $(\alpha ,\beta, 
\omega)$. $W^{cs}_n$ is $C^1$ in $\epsilon$ at point in the subset
$W^{cs}_{n+4}$. From Equation~\eqref{cspe1}, Claim~3 in the Theorem 
immediately follows.\hfill$\Box$

\begin{remark}\label{reop} Let $S^t$ denote the evolution 
operator of the perturbed nonlinear Schr\"odinger equation
\eqref{pnls}. The proofs of the Unstable Fiber Theorem and the 
Center-Stable Manifold Theorem also imply the following: $S^t$ is a 
$C^1$
map on $H^n$ for any fixed $t>0$, $n\geq 1$. $S^t$ is also $C^1$ in 
$(\alpha ,\beta ,\omega)$. $S^t$ is $C^1$ in $\epsilon$ as a map
from $H^{n+4}$ to $H^n$ for any fixed $n\geq 1$, $\epsilon \in 
[0,\epsilon_0)$, $\epsilon_0>0$.\end{remark}

\subsection{Stable Manifold of $Q_\epsilon$}

As mentioned earlier, the homoclinic orbit to be located will be 
asymptotic to the saddle $Q_\epsilon$ \eqref{Qec}. Dynamics on the
invariant plane $\Pi$ \eqref{Pi} on which $Q_\epsilon$ lives, is 
governed by Equations~\eqref{Ithe1}-\eqref{Ithe2} which are equivalent
to Equation~\eqref{nc1}-\eqref{nc2} with $f=0$. The eigenvalues of 
$Q_\epsilon$ are given by \eqref{Qee} on $\Pi$ and \eqref{leev} off
$\Pi$. Thus $Q_\epsilon$ has three unstable eigenvalues of two 
scales: One unstable eigenvalue of order 
$\mathcal{O}(\sqrt{\epsilon})$
with eigen-direction in $\Pi$, the other two unstable eigenvalues of 
order $\mathcal{O}(1)$ with eigen-directions off $\Pi$. On $\Pi$,
$Q_\epsilon$ has the unstable curve $\phi^u_{\sqrt{\e}}$ with approximate 
representation \eqref{cur}. Thus the $3D$ unstable manifold of
$Q_\epsilon $, $W^u(Q_\epsilon)$ has the representation
\begin{equation*} W^u(Q_\epsilon )=\cup_{p\in 
\phi^u_{\sqrt{\e}}}\mathcal{F}^+_p\end{equation*}
where $\mathcal{F}_p^+$ is the unstable fiber given in 
Theorem~\ref{UFT}. The scales of the stable eigenvalues of 
$Q_\epsilon$ range from
$\mathcal{O}(\epsilon)$ to $\mathcal{O}(\infty)$. The stable 
eigenvalue with eigen-direction in $\Pi$ has order
$\mathcal{O}(\sqrt{\epsilon})$. On $\Pi$, $Q_\epsilon$ has the stable 
curve $\phi^s_{\sqrt{\e}}$ with approximate representation \eqref{cur}.
 From the standard stable manifold theorem, $Q_\epsilon$ has a $C^1$ 
stable manifold $W^s_n(Q_\epsilon)$ in $H^n$ for any $n\geq 1$. In
fact, the codimension $3$ stable manifold of $Q_\epsilon$, 
$W^s_n(Q_\epsilon)$ intersects $\Pi$ along $\phi^s_{\sqrt{\e}}$. In order to 
locate a
homoclinic orbit, we need the size of $W^S_n(Q_\epsilon)$ large 
enough. Along $\phi^s_{\sqrt{\e}}$, the size of $W^s_n(Q_\epsilon)$ is
$\mathcal{O}(1)$ sufficient for our purpose. One can view 
$W^s_n(Q_\epsilon)$ as a wall with base $\phi^s_{\sqrt{\e}}$. As can be seen 
later in
the Second Measurement, one needs the size of $W^s_n(Q_\epsilon)$ off 
$\Pi$ to be of order $\mathcal{O}(\epsilon^\k)$, $\k <1$ in
order to overcome the order $\mathcal{O}(\epsilon)$ ``fuzz" between 
certain perturbed and unperturbed $(\epsilon =0)$ orbits to locate a
perturbed homoclinic orbit. Starting from the system 
\eqref{nc1}-\eqref{nc3}, one can only get the size of 
$W^s_n(Q_\epsilon)$ off $\Pi$
to be $\mathcal{O}(\epsilon)$ from standard stable manifold theorems. 
As discussed previously in the subsection on Normal Form
Transformation, an estimate of order $\mathcal{O}(\epsilon^\k)$, 
$\k<1$ can be achieved if the quadratic term $\mathcal{N}_2$
\eqref{wnc5} in \eqref{nc3} can be removed through a normal form 
transformation. Such a normal form transformation has be found in that
subsection.

\begin{theorem} The size of $W^s_n(Q_\epsilon)$ off $\Pi$ is of order 
$\mathcal{O}(\sqrt{\epsilon})$ for $\omega \in
\left(\frac{1}{2},\frac{3}{2}\right)/S$, where $S$ is a finite 
subset.\end{theorem}

\begin{proof} For $\omega \in \left( 
\frac{1}{2},\frac{3}{2}\right)/S$, where $S$ is a finite subset,
we apply the normal form transform given by 
\eqref{nors1}-\eqref{nors5} to Equation~\eqref{nc3}, then the system 
\eqref{nc1}-\eqref{nc3}
is transformed into the system \eqref{nfe1}-\eqref{nfe3}. By virtue 
of the estimate \eqref{wnfe}, the theorem follows from standard
argument. For details, see \cite{LMSW96}.
\end{proof}

As discussed in the subsection on Center-Stable Manifold, the 
center-stable manifold $W^{cs}_n$ is unique. Thus $W^s_n(Q_\epsilon)$ 
is a
codimension 1 submanifold of $W^{cs}_n$.

%% file: global.tex
\eqnsection{Global Theory}

Global Theory is referred to a theory global in phase space, which 
includes integrable theory, Melnikov measurement, and the so called
second measurements. These are tools necessary in locating a homoclinic orbit.

The entire process of locating the homoclinic orbit can be briefly 
summarized as follows: The integrable theory will provide explicit
representations for certain family of homoclinic orbits asymptotic to 
periodic orbits on the invariant plane $\Pi$. Local unstable fiber
theorem will provide ways of picking orbits in the local unstable 
manifold of $Q_\epsilon$, that are close to certain unperturbed
homoclinic orbits. Our main strategy is to use the unperturbed 
homoclinic orbits to trace the candidates for a perturbed homoclinic
orbit. The procedure is splinted into two steps:
\begin{enumerate}\item [Step 1.] Find an orbit that is in 
$W^u(Q_\epsilon)\cap W^{cs}_n$.
\item [Step 2.] Find out when this orbit is also in 
$W^u(Q_\epsilon)\cap W^s_n(Q_\epsilon)$, where $W^s_n(Q_\epsilon)$ is 
a codimension 1
submanifold of $W^{cs}_n$.\end{enumerate}
Step 1 will be accomplished through Melnikov measurement. The 
Melnikov vectors will be provided by integrable theory. The Melnikov
integrals will be evaluated along the unperturbed homoclinic orbits 
mentioned above. In contrast to the work \cite{LMSW96}, the new
feature in Step 1 is that $W^{cs}_n$ is not $C^1$ in $\epsilon$ 
everywhere rather only at its subset $W^{cs}_{n+4}$. This difficulty 
is
overcome by the fact that $W^u(Q_\epsilon)\subset H^n$ for any fixed 
$n\geq 1$ by virtue of the unstable fiber theorem. Step 2 will be
accomplished by the so called second measurement. It turns out that 
one can trace the perturbed orbit in $W^u(Q_\epsilon)\cap W^{cs}_n$
through an unperturbed homoclinic orbit to an order 
$\mathcal{O}(\epsilon |\ln \epsilon |)$ neighborhood of $\Pi$ 
\eqref{Pi}.
In order to check when this orbit can be in $W^u(Q_\epsilon)\cap 
W^s_n(Q_\epsilon )$, one needs the size of $W^s_n(Q_\epsilon)$ off 
$\Pi$
to be large enough, and $\mathcal{O}(\sqrt{\epsilon})$ is sufficient.

\subsection{Integrable Theory}

Consider the integrable 1D cubic focusing nonlinear Schr\"odinger 
equation ($\e = 0$ in (\ref{pnls})),
\begin{equation}
iq_t = q_{xx} + 2[|q|^2 - \om^2] q \ .
\label{NLS}
\end{equation}
Its Lax pair is given by the Zakharov-Shabat linear system,
\begin{eqnarray}
\psi_x &=& U \psi \ , \label{ZS1} \\
\psi_t &=& V \psi \ , \label{ZS2}
\end{eqnarray}
where
\[
U = i \left ( \begin{array}{lr} \la & q \cr \bq & -\la \cr 
\end{array} \right ) \ , 
\]
\[
V = i \left ( \begin{array}{lr} 2\la^2 -|q|^2 +\om^2 & 
2\la q -i q_x \cr 2 \la \bq + i \overline{q_x} & -2 \la^2 +|q|^2 -\om^2 \cr 
\end{array} \right ) \ . 
\]

\subsubsection{Isospectral Theory}

Focusing one's attention on the spatial part (\ref{ZS1}) of the Lax 
pair (\ref{ZS1},\ref{ZS2}), one can define the fundamental matrix 
solution $M(x)$, s.t. $M(0)$ is the $2\times 2$ identity matrix. Then 
the Floquet discriminant $\Dl$ is defined as
\[
\Dl = \ \mbox{trace}\ M(2\pi)\ .
\]
$\Dl = \Dl (\la, q)$, as a functional in $q$ for any $\la \in \mathbb{C}$, 
provides enough functionally independent constants of motion to make 
NLS (\ref{NLS}) integrable in the classical Liouville sense. For each 
fixed $q$, there is a sequence of special points 
$\{ \la^s_j, \ j \in Z \}$ of $\la \in \mathbb{C}$ called simple points for 
which $|\Dl(\la^s_j,q)| = 2$. There is also a sequence of critical 
points $\{ \la^c_j, \ j \in Z \}$ of $\la \in \mathbb{C}$ for 
which $\frac {\pa}{\pa \la}\Dl(\la^c_j,q) = 0$. When some $\la^s_j$ 
coincides with some $\la^c_l$, a double point is formed. The geometric 
multiplicity is the dimension of the eigenspace of (\ref{ZS1}) at the 
double point. 
\begin{definition}
The sequence of constants of motion $F_j$ is defined as
\begin{equation}
F_j = \Dl(\la^c_j,q)\ , \ \ j \in Z\ .
\label{Fj}
\end{equation}
\end{definition}
\nid
$F_j$'s provide a sequence of Melnikov functions. More importantly,
the gradients of $F_j$'s, which will be the Melnikov vectors,  have a simple representation,
\begin{equation}
\frac{\dl F_j}{\dl \vq} = i \frac {\sqrt{\Dl^2 -4}} {W(\psi^+, \psi^-)}
\left ( \begin{array}{c} \psi_2^+ \psi_2^- \cr - \psi_1^+ \psi_1^-
\cr \end{array} \right ) \ , \ 
\ \mbox{at}\ \la = \la^c_j, \ \ j \in Z\ ,
\label{gFj}
\end{equation}
where $\vq = (q,\bq)^T$, $\psi^{\pm}=(\psi^{\pm}_1, \psi^{\pm}_2)^T$ 
are two eigenfunctions at $\la = \la^c_j$, and $W(\psi^+,\psi^-)=
\psi^+_1 \psi^-_2 - \psi^+_2 \psi^-_1$ is the Wronskian.
For more details on the isospectral theory of NLS, we refer the readers to \cite{LM94}.

\subsubsection{B\"acklund-Darboux Transformation}

The particular form of the B\"acklund-Darboux transformation for 
NLS (\ref{NLS}), that is useful for our purpose, is due to David 
Sattinger and V. Zurkowski \cite{SZ87}.
\begin{theorem}
Let $q(t,x)$ be a solution of NLS (\ref{NLS}), $\nu$ is a 
complex double point of geometric multiplicity $2$. Let $\phi^{\pm}$
be two linearly independent eigenfunctions of the Lax pair (\ref{ZS1},\ref{ZS2}) at $\la = \nu$. Denote by $\phi$ the general 
solution
\[
\phi = \phi(t, x, \nu, c_+, c_-) = c_+ \phi^+ + c_- \phi^-\ ,
\]
We use $\phi$ to define a Gauge transformation matrix 
\begin{equation}
G = G(\la; \nu, \phi) = \Ga \left ( \begin{array}{lr} \la -\nu & 0 \cr 
0 & \la - \bar{\nu} \cr \end{array} \right )  \Ga^{-1} \ ,
\label{Gm}
\end{equation}
where 
\[
\Ga = \left ( \begin{array}{lr} \phi_1 & - \overline{\phi_2} \cr 
\phi_2 & \overline{\phi_1} \cr  \end{array} \right ) \ .
\]
Then we define $Q$ and $\Psi$ by
\[
Q = q + 2 (\nu - \bar{\nu}) \frac {\phi_1 \overline{\phi_2}} {|\phi_1|^2 
+ |\phi_2|^2} \ ,
\]
and 
\[
\Psi = G \psi \ ,
\]
where $\psi$ solves the Lax pair (\ref{ZS1},\ref{ZS2}) at ($\la, q$).
Then $\Psi$ solves the Lax pair (\ref{ZS1},\ref{ZS2}) at ($\la, Q$),
and $Q$ also solves NLS (\ref{NLS}).
\label{bdthm}
\end{theorem}
                                
\subsubsection{Figure Eight Structures}

Consider the special solution of NLS (\ref{NLS}),
\begin{equation}
q_c = a e^{i\th(t)}\ , \ \ \th(t) = -[2(a^2 - \om^2)t +\ga ]\ .
\label{usoln}
\end{equation}
The corresponding Floquet discriminant is given by
\[
\Dl(\la, q_c) = 2 \cos (2\pi k) \ , \ \ k = \sqrt{a^2 + \la^2}\ ,
\]
and two eigenfunctions (Bloch functions) are 
\begin{equation}
\psi^{\pm} = \left ( \begin{array}{c} a e^{i \frac {\th}{2}} \cr 
(\pm k - \la)e^{-i \frac {\th}{2}} \cr \end{array} \right ) 
\exp \{ \pm i2 \la k t \pm ikx \} \ .
\label{Bf}
\end{equation}
When $k$ is real, to have temporal growth (and decay) in $\psi^{\pm}$, one needs $\la$ to be purely imaginary. The temporal growth (and decay) in $\psi^{\pm}$ is connected to the linear instability of $q_c$, since 
quadratic products of $\psi^{\pm}$ solve linearized NLS \cite{LM94}. 
The temporal growth is also necessary for constructing homoclinic 
solutions through the B\"acklund-Darboux transformation. Specifically,
the double points of $\Dl$ are given by
\[
k = \sqrt{a^2 + \la^2} = j/2\ , \ \ j \in Z/\{0\} \ .
\]
If one requires that $a$ lies in the interval
\[
a \in (1/2, 1)\ , 
\]
then there is only one pair of complex double points
\[
\la = \pm \nu = \pm i \sg \ , \ \ \sg = \sqrt{a^2 - 1/4}\ .
\]
If one requires that $a$ lies in the interval
\[
a \in (1, 3/2)\ , 
\]
then there are two pairs of complex double points
\[
\la = \pm \nu = \pm i \sg \ , \ \mbox{and} \ \la = \pm \hat{\nu} 
= \pm i \hat{\sg}\ , \hat{\sg} = \sqrt{a^2 - 1}\ .
\]
Next we will construct homoclinic orbits, starting from the special 
solution $q_c$, through the B\"acklund-Darboux transformation. Notice that building the B\"acklund-Darboux transformation at $\la = \nu$ v.s. at $\la = - \nu$ and 
at $\la = \hat{\nu}$ v.s. at $\la = - \hat{\nu}$ always 
lead to equivalent results. We will choose $\la = \nu$ and 
$\la = \hat{\nu}$. \\

\nid
{\bf One Pair of Complex Double Points Case} \\

Let $\phi^{\pm} = \psi^{\pm}(t, x, \nu)$ defined in (\ref{Bf}), and let 
\begin{equation}
\phi = c^+ \phi^+ + c^- \phi^- \ .
\label{Gf}
\end{equation}
Applying the B\"acklund-Darboux transformation given in Theorem \ref{bdthm}, one gets a new solution,
\begin{eqnarray}
Q &=& q_c \bigg [ 1 + \sin \vth_0 \ \mbox{sech} \tau 
\cos y \bigg ]^{-1} \cdot \bigg [ \cos 2\vth_0 - i \sin 2\vth_0 \tanh \tau
 \nonumber \\
& & - \sin \vth_0 \ \mbox{sech} \tau \cos y \bigg ] \ , \label{sne}
\end{eqnarray}
where
\begin{equation}
c^+/c^- = e^{\rho +i \vth}\ , \ \ \frac{1}{2} + \nu = a e^{i\vth_0}\ ,
\ \ \tau = 2 \sg t - \rho \ , \ \ y = x + \vth - \vth_0 +\pi/2\ .
\label{par1}
\end{equation}
As $t \ra \pm \infty$, 
\begin{equation}
Q \ra q_c e^{\mp i2\vth_0}\ .
\label{asy1}
\end{equation}
Thus $Q$ is asymptotic to $q_c$ up to phase shifts as $t \ra \pm \infty$.
We say $Q$ is a homoclinic orbit asymptotic to the periodic orbit given 
by $q_c$. For a fixed amplitude $a$ of $q_c$, the phase $\ga$ of $q_c$ and 
the B\"acklund parameters $\rho$ and $\vth$ parametrize a $3$-dimensional 
submanifold with a figure eight structure. For an illustration, see 
Figure \ref{snef}. 
\begin{figure}
\vspace{1.5in}
\caption{Figure eight structure of noneven data with one unstable mode.}
\label{snef}
\end{figure}
\begin{figure}
\vspace{1.5in}
\caption{Figure eight structure of even data with one unstable mode.}
\label{sef}
\end{figure}
\nid 
If one restricts the 
B\"acklund parameter $\vth$ by $\vth -\vth_0 +\pi/2 = 0$, or $\pi$, 
one gets $Q$ to be even in $x$, 
\begin{eqnarray}
Q &=& q_c \bigg [ 1 \pm \sin \vth_0 \ \mbox{sech} \tau  
\cos x \bigg ]^{-1} \nonumber \\
& & \cdot \bigg [ \cos 2\vth_0 - i \sin 2\vth_0 \tanh \tau \mp 
\sin \vth_0 \ \mbox{sech} \tau \cos x \bigg ] \ , \label{se}
\end{eqnarray}
where the upper sign corresponds to $0$. Then for a fixed amplitude 
$a$ of $q_c$, the phase $\ga$ of $q_c$ and 
the B\"acklund parameter $\rho$ parametrize a $2$-dimensional 
submanifold with a figure eight structure. For an illustration, see 
Figure \ref{sef}. \\

\nid
{\bf Two Pairs of Complex Double Points Case} \\

Let $\hat{\phi}^{\pm} = \psi^{\pm}(t, x, \hat{\nu})$ defined in (\ref{Bf}), and let 
\begin{equation}
\hat{\phi} = \hat{c}^+ \hat{\phi}^+ + \hat{c}^- \hat{\phi}^- \ .
\label{hGf}
\end{equation}
In this ``two pairs of complex double points'' case, to get the complete 
foliation of the figure eight structure, one needs to iterate the 
B\"acklund-Darboux transformation. First one needs to apply the 
B\"acklund-Darboux transformation at $\la = \nu$, then one needs to 
iterate the B\"acklund-Darboux transformation at $\la = \hat{\nu}$.
Switching the order between $\nu$ and $\hat{\nu}$ leads to the same 
result. At $\la = \nu$, the Gauge transform $G = G(\la; \nu, \phi)$ 
(\ref{Gm}), then one defines
\begin{equation}
\hat{\Phi}^{\pm} = G(\hat{\nu}; \nu, \phi)\hat{\phi}^{\pm}\ .
\label{hPhipm}
\end{equation}
Let 
\begin{equation}
\hat{\Phi} = G(\hat{\nu}; \nu, \phi)\hat{\phi} = \hat{c}^+ \hat{\Phi}^+ + 
\hat{c}^- \hat{\Phi}^- \ .
\label{iGf}
\end{equation}
After an iteration on the B\"acklund-Darboux transformation, one gets 
the solution of NLS (\ref{NLS}) with the representation,
\[
\tilde{Q} = q_c +2(\nu - \bar{\nu}) \frac {\phi_1 \overline{\phi_2}} 
{|\phi_1|^2+|\phi_2|^2} + 2(\hat{\nu} - \bar{\hat{\nu}}) \frac {\hat{\Phi}_1 
\overline{\hat{\Phi}_2}} {|\hat{\Phi}_1|^2+|\hat{\Phi}_2|^2}\ .
\]
Explicit formula for $\tilde{Q}$ is,
\begin{equation}
\tilde{Q} = Q + q_c \frac{\W_2 \sin \hvth_0}{\W_1}\ ,
\label{dne}
\end{equation}
where $Q$ is given in (\ref{sne}),
\begin{eqnarray*}
\W_1 &=& \bigg [ (\sin \hvth_0)^2(1+\sin \vth_0 \ \mbox{sech} \tau 
\cos y )^2 +\frac{1}{8} (\sin 2\vth_0)^2 (\mbox{sech} \tau)^2 
(1 -\cos 2y) \bigg ] \\
& & \cdot (1 + \sin \hvth_0 \ \mbox{sech} \htau 
\cos \hy ) \\
&-& \frac {1}{2} \sin 2\vth_0 \sin 2\hvth_0 \ \mbox{sech} \tau 
\ \mbox{sech} \htau (1+\sin \vth_0 \ \mbox{sech} \tau 
\cos y ) \sin y \sin \hy  \\
&+& (\sin \vth_0)^2 \bigg [ 1+ 2 \sin \vth_0 \ \mbox{sech} \tau 
\cos y + [(\cos y)^2 - (\cos \vth_0)^2](\mbox{sech} \tau)^2 \bigg ]
\\
& &\cdot (1 + \sin \hvth_0 \ \mbox{sech} \htau 
\cos \hy ) \\
&-& 2\sin \hvth_0 \sin \vth_0 \bigg [ \cos \hvth_0 \cos \vth_0
\tanh \htau \tanh \tau + ( \sin \vth_0 + \ \mbox{sech} \tau 
\cos y)\\
& & \cdot ( \sin \hvth_0 + \ \mbox{sech} \htau 
\cos \hy) \bigg ] (1 + \sin \vth_0 \ \mbox{sech} \tau 
\cos y ) \ ,
\end{eqnarray*}
\begin{eqnarray*}
\W_2 &=& \bigg [ -2 (\sin \hvth_0)^2 (1 + \sin \vth_0 \ 
\mbox{sech} \tau \cos y )^2 +\frac {1}{4} (\sin 2\vth_0)^2 
(\mbox{sech} \tau)^2 (1-\cos 2y)\bigg ]\\
& & \cdot (\sin \hvth_0 + 
\ \mbox{sech} \htau \cos \hy + i \cos \hvth_0 \tanh \htau ) \\
&+& 2 (\sin \vth_0)^2(-\cos \vth_0 \tanh \tau + i \sin \vth_0 + 
i\ \mbox{sech} \tau \cos y)^2 \\
& & \cdot (\sin \hvth_0 + 
\ \mbox{sech} \htau \cos \hy - i \cos \hvth_0 \tanh \htau )  \\
&+& 2 \sin \vth_0 (\sin \vth_0 + 
\ \mbox{sech} \tau \cos y + i \cos \vth_0 \tanh \tau ) 
\\
& & \cdot \bigg [2 \sin \hvth_0 (1 + \sin \vth_0 \ 
\mbox{sech} \tau \cos y )(1 + \sin \hvth_0 \ 
\mbox{sech} \htau \cos \hy )\\
& & - \sin 2\vth_0 \cos \hvth_0 
\ \mbox{sech} \tau \ \mbox{sech} \htau  \sin y \sin \hy \bigg ] \ ,
\end{eqnarray*}
and the notations are given by 
\[
1 + \hat{\nu} = a e^{i\hvth_0} \ , \ \ \hat{c}^+/\hat{c}^- 
=e^{\hat{\rho}+i \hvth}\ , \ \ \hat{\tau} = 4 \hat{\sg} t - \hat{\rho} \ ,
\ \ \hy = 2x + \hvth - \hvth_0 +\pi/2\ .
\]
The asymptotic phase of $Q$ is as follows, as $t \ra \pm \infty$,
\begin{equation}
\tilde{Q} \ra q_c e^{\mp i 2 (\vth_0 + \hvth_0)}\ . 
\label{dasym}
\end{equation}
Thus $\tilde{Q}$ is asymptotic to $q_c$ up to phase shifts as 
$t \ra \pm \infty$. We say $\tilde{Q}$ is a homoclinic orbit 
asymptotic to the periodic orbit given by $q_c$. For a fixed amplitude 
$a$ of $q_c$, the phase $\ga$ of $q_c$ and the B\"acklund parameters $\rho$, $\vth$, $\hat{\rho}$, and $\hvth$ parametrize a $5$-dimensional 
submanifold with a figure eight structure. For an illustration, see 
Figure \ref{dnef}. 
\begin{figure}
\vspace{1.5in}
\caption{Figure eight structure of noneven data with two unstable modes.}
\label{dnef}
\end{figure}
\begin{figure}
\vspace{1.5in}
\caption{Figure eight structure of even data with two unstable modes.}
\label{def}
\end{figure}
If one put restrictions on the B\"acklund parameters $\vth$ and 
$\hvth$, s.t. 
\begin{equation}
\vth - \vth_0 +\pi/2 = \left \{ \begin{array}{c} 0 
\ \left \{ \begin{array}{c} \hvth - \hvth_0 +\pi/2 = 0 \ ,\cr
\hvth - \hvth_0 +\pi/2 = \pi \ , \cr \end{array} \right.
\cr \pi \ \left \{ \begin{array}{c} \hvth - \hvth_0 +\pi/2 = 0 \ ,\cr
\hvth - \hvth_0 +\pi/2 = \pi \ , \cr \end{array} \right.
\cr \end{array} \right.
\label{dec}
\end{equation}
then $Q$ is even in $x$. Then for a fixed amplitude 
$a$ of $q_c$, the phase $\ga$ of $q_c$ and the B\"acklund parameters $\rho$ and $\hat{\rho}$ parametrize a $3$-dimensional 
submanifold with a figure eight structure. For an illustration, see 
Figure \ref{def}. 

\subsubsection{Melnikov Vectors}

Notice that (\ref{gFj}) evaluated at $(\nu, Q)$ and $(\bar{\nu}, Q)$
are linearly dependent. Same is true for $(\nu, \tilde{Q})$ or 
$(\hat{\nu}, \tilde{Q})$. \\

\nid
{\bf One Pair of Complex Double Points Case} \\

In this case, the Melnikov vector is $\frac{\dl F_1}{\dl \vq}$, (\ref{gFj}) 
at $\la = \nu$, evaluated along the homoclinic orbit $Q$ (\ref{sne}) 
or (\ref{se}). 
\[
\frac{\dl F_1}{\dl \vq} = i \frac {\sqrt{\Dl^2(\nu) -4}} {W(\Phi^+, \Phi^-)}
\left ( \begin{array}{c} \Phi_2^+ \Phi_2^- \cr - \Phi_1^+ \Phi_1^-
\cr \end{array} \right ) \ , 
\]
where (cf: (\ref{Gm})),
\begin{eqnarray}
\Phi^{\pm} &=& G(\nu; \nu, \phi)\phi^{\pm}\label{Phipm}\\
&=& \pm c^{\mp} W(\phi^+, \phi^-)\frac {\nu - \bar{\nu}} {|\phi_1|^2 + 
|\phi_2|^2} \left (\begin{array}{c} \overline{\phi_2} \cr -\overline{\phi_1} 
\cr \end{array} \right ) \ . \nonumber
\end{eqnarray}
By L'Hospital's rule,
\[
\frac {\sqrt{\Dl^2 -4}} {W(\Phi^+, \Phi^-)} = 
\frac {\sqrt{\Dl(\nu)\Dl''(\nu)}} {(\nu - \bar{\nu})W(\phi^+, \phi^-)}\ .
\]
$\phi$ (\ref{Gf}) can be rewritten as
\begin{equation}
\phi_1 = 2\sqrt{c^+c^-} a e^{i\th/2} u_1\ , \ \ 
\phi_2 = 2\sqrt{c^+c^-} a e^{-i\th/2} u_2\ , 
\label{rew1}
\end{equation}
where
\begin{eqnarray*}
u_1 &=&  \cosh \frac{\tau}{2} \cos z -i \sinh \frac{\tau}{2} \sin z \ , \\
u_2 &=&  -\sinh \frac{\tau}{2} \cos (z - \vth_0) + i \cosh  \frac{\tau}{2}
\sin (z - \vth_0) \ ,
\end{eqnarray*}
where
\[
z = x/2 + \vth/2 \ ,
\]
and other notations have been defined in (\ref{usoln}, \ref{par1}, 
\ref{dne}). Finally, one gets the explicit representation for 
the Melnikov vector,
\begin{equation}
\frac{\dl F_1}{\dl \vq} = \frac{1}{4} a^{-2} i (\nu - \bar{\nu})
\sqrt{\Dl(\nu)\Dl''(\nu)} \frac {1}{(|u_1|^2 + |u_2|^2)^2}
\left ( \begin{array}{c} \overline{q_c}\ \overline{u_1}^{\ 2} \cr 
- q_c\ \overline{u_2}^{\ 2} \cr \end{array} \right ) \ .
\end{equation} \\

\nid
{\bf Two Pairs of Complex Double Points Case} \\

In this case, the Melnikov vectors are $\frac{\dl F_1}{\dl \vq}$ and $\frac{\dl F_2}{\dl \vq}$ (\ref{gFj}) at $\la = \nu$ and 
$\la = \hat{\nu}$ respectively, evaluated along the homoclinic 
orbit $\tilde{Q}$ (\ref{dne}) or (\ref{dec}). We know that 
$\hat{\Phi}$ is defined in (\ref{iGf}). Then we use $\hat{\Phi}$ 
to define a Gauge matrix $G(\la; \hat{\nu}, \hat{\Phi})$. Let
\begin{equation}
\Phi^{(\pm,*)} = G(\nu; \hat{\nu}, \hat{\Phi}) \Phi^{\pm} \ , 
\ \ \hat{\Phi}^{(\pm,*)} = G(\hat{\nu}; \hat{\nu}, \hat{\Phi}) 
\hat{\Phi}^{\pm} \ , 
\label{stphi}
\end{equation}
where $\Phi^{\pm}$ and $\hat{\Phi}^{\pm}$ are defined in 
(\ref{Phipm}, \ref{hPhipm}). Then the Melnikov vectors are
\begin{eqnarray}
\frac{\dl F_1}{\dl \vq} &=& i \frac {\sqrt{\Dl^2(\nu) -4}} {W(\Phi^{(+,*)}, \Phi^{(-,*)})}
\left ( \begin{array}{c} \Phi_2^{(+,*)}\Phi_2^{(-,*)} \cr 
- \Phi_1^{(+,*)}\Phi_1^{(-,*)}\cr \end{array} \right ) \ , 
\label{gF1} \\
\frac{\dl F_2}{\dl \vq} &=& i \frac {\sqrt{\Dl^2(\hat{\nu}) -4}} {W(\hat{\Phi}^{(+,*)}, \hat{\Phi}^{(-,*)})}
\left ( \begin{array}{c} \hat{\Phi}_2^{(+,*)}\hat{\Phi}_2^{(-,*)} \cr 
- \hat{\Phi}_1^{(+,*)}\hat{\Phi}_1^{(-,*)}\cr \end{array} \right ) \ , 
\label{gF2}
\end{eqnarray}
By L'Hospital's rule,
\[
\frac {\sqrt{\Dl^2(\nu) -4}} {W(\Phi^{(+,*)}, \Phi^{(-,*)})}=
\frac {\sqrt{\Dl(\nu)\Dl''(\nu)}} {(\nu - \bar{\nu})
(\nu - \hat{\nu})(\nu - \bar{\hat{\nu}})W(\phi^+, \phi^-)}\ ,
\]
\[
\frac {\sqrt{\Dl^2(\hat{\nu}) -4}} 
{W(\hat{\Phi}^{(+,*)}, \hat{\Phi}^{(-,*)})}=
\frac {\sqrt{\Dl(\hat{\nu})\Dl''(\hat{\nu})}} {(\hat{\nu} - \bar{\hat{\nu}})(\hat{\nu}- \nu )(\hat{\nu} - \bar{\nu})
W(\hat{\phi}^+, \hat{\phi}^-)}\ .
\]
We know that $\phi$ (\ref{Gf}) can be rewritten as (\ref{rew1}).
$\hat{\phi}$ (\ref{hGf}) can also be rewritten as
\begin{equation}
\hat{\phi}_1 = 2\sqrt{\hat{c}^+\hat{c}^-} a e^{i\th/2} v_1\ , \ \ 
\hat{\phi}_2 = 2\sqrt{\hat{c}^+\hat{c}^-} a e^{-i\th/2}v_2\ , 
\label{rew2}
\end{equation}
where
\begin{eqnarray*}
v_1 &=&  \cosh \frac{\hat{\tau}}{2} \cos \hat{z} -i \sinh \frac{\hat{\tau}}{2} \sin \hat{z} \ , \\
v_2 &=&  -\sinh \frac{\hat{\tau}}{2} \cos (\hat{z} - \hvth_0) + i 
\cosh \frac{\hat{\tau}}{2} \sin (\hat{z}- \hvth_0) \ ,
\end{eqnarray*}
where
\[
\hat{z} = x + \hvth/2 \ ,
\]
and other notations have been defined in (\ref{usoln}, \ref{par1}, 
\ref{dne}). Using (\ref{rew1}, \ref{rew2}), one can get the representation for $\hat{\Phi}$ (\ref{iGf}),
\[
\hat{\Phi}_1 = 2\sqrt{\hat{c}^+\hat{c}^-} a e^{i\th/2} V_1\ , \ \ 
\hat{\Phi}_2 = 2\sqrt{\hat{c}^+\hat{c}^-} a e^{-i\th/2}V_2\ , 
\]
where $V_1$ and $V_2$ are defined as
\begin{eqnarray}
V_1 &=& \frac {1} {|u_1|^2+|u_2|^2} \bigg [ [(\hat{\nu}- \nu )
|u_1|^2 + (\hat{\nu}- \bar{\nu} )|u_2|^2] v_1 + (\bar{\nu}- \nu )
u_1 \overline{u_2} v_2 \bigg ]\ , \label{V1} \\
V_2 &=& \frac {1} {|u_1|^2+|u_2|^2} \bigg [ (\bar{\nu}- \nu )
\overline{u_1} u_2 v_1 + [(\hat{\nu}- \bar{\nu})
|u_1|^2 + (\hat{\nu}- \nu )|u_2|^2] v_2 \bigg ]\ . \label{V2}
\end{eqnarray}
Finally, one gets the explicit representations
\begin{eqnarray}
\frac{\dl F_1}{\dl \vq} &=& \frac{1}{4} a^{-2} i (\nu - \bar{\nu})  
(\nu - \hat{\nu})^{-1}(\nu - \bar{\hat{\nu}})^{-1}
\sqrt{\Dl(\nu)\Dl''(\nu)} \left ( \begin{array}{c} \overline{q_c} S^2_2 \cr 
- q_c S^2_1 \cr \end{array} \right )\ , \label{exF1} \\
\frac{\dl F_2}{\dl \vq} &=& \frac{1}{2} a^{-2} i (\hat{\nu} - \bar{\hat{\nu}})(\hat{\nu}- \nu )(\hat{\nu} - \bar{\nu})
\sqrt{\Dl(\hat{\nu})\Dl''(\hat{\nu})} \left ( \begin{array}{c} \overline{q_c} \hat{S}^2_2 \cr 
- q_c \hat{S}^2_1 \cr \end{array} \right )\ , \label{exF2}
\end{eqnarray}
where $S_l$ and $\hat{S}_l$ ($l=1,2$) are independent of the phase 
$\ga$ of $q_c$, and have the representations
\begin{eqnarray}
S_1 &=& \frac {1} {(|u_1|^2+|u_2|^2)(|V_1|^2+|V_2|^2)} 
\bigg [ [(\nu - \hat{\nu})
|V_1|^2 + (\nu - \bar{\hat{\nu}} )|V_2|^2] \overline{u_2} \nonumber \\
& & - (\bar{\hat{\nu}}- \hat{\nu} )
V_1 \overline{V_2} \overline{u_1}\bigg ]\ , \label{S1} \\
S_2 &=& \frac {1} {(|u_1|^2+|u_2|^2)(|V_1|^2+|V_2|^2)} 
\bigg [ (\bar{\hat{\nu}}- \hat{\nu} ) \overline{V_1} V_2 \overline{u_2} -
[(\nu - \bar{\hat{\nu}})|V_1|^2 \nonumber \\
& & + (\nu - \hat{\nu} )|V_2|^2] 
\overline{u_1} \bigg ]\ , \label{S2} \\
\hat{S}_1 &=& \frac {\overline{V_2}} {|V_1|^2+|V_2|^2}\ , \label{hS1} \\
\hat{S}_2 &=& \frac {\overline{V_1}} {|V_1|^2+|V_2|^2}\ . \label{hS2} 
\end{eqnarray}

\subsection{Melnikov Analysis}

Let $p$ be any point on $\phi^u_{\sqrt{\e}}$ \eqref{cur} which is the 
unstable curve of $Q_\epsilon $ in $\Pi$ (\ref{Pi}). Let $q_\epsilon (0)$ 
and $q_0(0)$ be any two points on the unstable fibers 
$\mathcal{F}_p^+\mid_\epsilon$ and 
$\mathcal{F}_p^+\mid_{\epsilon=0}$, with the same $\xi^+_k$
coordinates. By the Unstable Fiber Theorem, $\mathcal{F}^+_p$ is 
$C^1$ in $\epsilon $ for $\epsilon \in [0,\epsilon_0)$, $\epsilon_0
>0$, thus
\begin{equation*}\| q_\epsilon (0)-q_0(0)\|_{n+8}\leq C\epsilon.\end{equation*}
The key point here is that $\mathcal{F}^+_p\subset H^s$ for any fixed 
$s\geq 1$. By Remark~\ref{reop}, the evolution operator of the
perturbed NLS equation \eqref{pnls} $S^t$ is $C^1$ in $\epsilon$ as a 
map from $H^{n+4}$ to $H^n$ for any fixed $n\geq 1$, $\epsilon \in
[0,\epsilon_0)$, $\epsilon _0>0$. Also $S^t$ is a $C^1$ map on $H^n$ 
for any fixed $t>0$, $n\geq 1$. Thus
\begin{equation*}\| q_\epsilon (T)-q_0(T)\| _{n+4}=\| S^T(q_\epsilon 
(0))-S^T(q_0(0))\|_{n+4}\leq C_1\epsilon,\end{equation*}
where $T>0$ is large enough so that
\begin{equation*}q_0(T)\in W^{cs}_{n+4}\mid_{\epsilon=0}.\end{equation*}
Our goal is to determine when $q_\epsilon (T)\in W^{cs}_n$ through 
Melnikov measurement. Let $q_\epsilon (T)$ and $q_0(T)$ have the coordinate
expressions
\begin{equation}\label{sme} q_\epsilon 
(T)=(\xi^{+,\epsilon}_k,v_\epsilon),\quad 
q_0(T)=(\xi^{+,0}_k,v_0).\end{equation}
Let $\tilde{q}_\epsilon (T)$ be the unique point on 
$W^{cs}_{n+4}$, 
which has the same $v$-coordinate as $q_\epsilon (T)$,
\begin{equation*}\tilde{q}_\epsilon 
(T)=(\tilde{\xi}^{+,\epsilon}_k,v_\epsilon)\in 
W^{cs}_{n+4}.\end{equation*}
By the Center-Stable Manifold Theorem, at points in the subset 
$W^{cs}_{n+4}$, $W^{cs}_n$ is $C^1$ smooth in $\epsilon$ for $\epsilon
\in [0,\epsilon_0)$, $\epsilon_0 >0$, thus
\begin{equation}\label{sme1} \| q_\epsilon (T)-\tilde{q}_\epsilon 
(T)\|_n\leq C_2\epsilon .\end{equation}
Also our goal now is to determine when the signed distances
\begin{equation*}\xi^{+,\epsilon}_k-\tilde{\xi}^{+,\epsilon}_k,\quad 
(k=1,2)\end{equation*}
are zero through Melnikov measurement. Equivalently, one can define 
the signed distances
\begin{equation*}\begin{split} d_k&=\langle \nabla 
F_k(q_0(T)),q_\epsilon (T)-\tilde{q}_\epsilon (T)\rangle\\
&\equiv \partial_qF_k(q_0(T)) (q_\epsilon (T)-\tilde{q}_\epsilon (T))\\
&+ \partial_{\bar{q}}F_k(q_0(T)) (q_\epsilon 
(T)-\tilde{q}_\epsilon (T))^-,\quad k=1,2,\end{split}\end{equation*}
where $F_k$ and $\nabla F_k$ are given in the subsection on 
Integrable Theory, $q_0(t)$ is the homoclinic orbit also given in the 
same subsection. In fact, $q_\epsilon (t)$, $\tilde{q}_\epsilon (t)$, 
$q_0(t)\in H^n$, for any fixed $n\geq 1$. The rest of the derivation
for Melnikov integrals is completely standard. For details, see 
\cite{LMSW96} \cite{LM97}.
\begin{equation}\label{Meln1} d_k=\epsilon M_k+o(\epsilon ),\quad 
k=1,2,\end{equation}
where
\begin{equation*}\begin{split} 
M_k&=\int^{+\infty}_{-\infty}\int^{2\pi}_0 
[\partial_qF_k(q_0(t))(\partial^2_xq_0(t)-\alpha
q_0(t)+\beta)\\
&\quad + 
\partial_{\bar{q}}F_k(q_0(t))(\partial^2_x\overline{q_0(t)}-\alpha 
\overline{q_0(t)}+\beta )]dxdt,\end{split}\end{equation*}
where $q_0(t)$, $\partial_qF_k$, and $\partial_{\bar{q}}F_k$ are 
given in the subsection on Integrable Theory.

\begin{theorem}\label{Melthm} There exists $\epsilon_0>0$, such that 
for any $\e \in (0, \e_0)$, there exists a domain $\mathcal{D}_\e 
\subset \mathbb{R}^+ \times \mathbb{R}^+ \times \mathbb{R}^+$ where 
$\om \in (\frac{1}{2}, \frac{3}{2})/S$, $S$ is a finite subset, and 
$\al \om < \be$. For any $(\alpha ,\beta, \omega) \in \mathcal{D}_\e$,
there exists another orbit in $W^u(Q_\epsilon )\cap W^{cs}_n$ other 
than the unstable curve $\phi^u_{\sqrt{\e}}$ (\ref{cur}) of $Q_\e$, for 
the perturbed nonlinear Schr\"odinger equation \eqref{pnls}.
\end{theorem}
\begin{proof} This theorem follows immediately from the explicit 
computation in the subsection on Evaluation of Melnikov Integrals and 
Second Distance, and the implicit function theorem. 
\end{proof}

\subsection{The Second Measurement}

The second measurement starts with the orbit obtained in 
Theorem~\ref{Melthm}, i.e. $q_\epsilon(t)$ where
$q_\epsilon(T)=\tilde{q}_\epsilon (T)$. The goal is to determine when 
$q_\epsilon (t)$ is also in $W^s_n(Q_\epsilon)$. Recall that
$W^s_n(Q_\epsilon)$ can be visualized as a codimension-one wall in 
$W^{cs}_n$ with base curve in $\Pi$ and with 
$\mathcal{O}(\sqrt{\epsilon})$
height. Thus we have to continue to follow $q_\epsilon (t)$ and 
$q_0(t)$ to a smaller neighborhood of $\Pi$. From the explicit
expression of $q_0(t)$, we know that $q_0(t)$ approaches $\Pi$ at the 
rate $\mathcal{O}(e^{-\mu t})$,
\begin{equation}\mu=\min \{ 
\sqrt{4\omega^2-1},4\sqrt{\omega^2-1}\}\end{equation}
(cf: \eqref{leev}). Thus
\begin{equation}\label{sarg1} \text{distance} \left\{ 
q_0(T+\frac{1}{\mu}|\ln \epsilon |), \Pi\right\}<C\epsilon.\end{equation}

\begin{lemma}\label{lesm} For all $t\in \left[T, T+\frac{1}{\mu}|\ln 
\epsilon |\right]$,
\begin{equation}\| q_\epsilon (t)-q_0(t)\|_n\leq \tilde{C}_1\epsilon |\ln 
\epsilon |^2,\end{equation}
where $\tilde{C}_1 =\tilde{C}_1(T)$.
 \end{lemma}

\begin{proof} We start with the system \eqref{cspe1}-\eqref{cspe2}. Let
\begin{align*} q_\epsilon (t)&= (\xi^{+,\epsilon}_k(t),J^\epsilon 
(t),\theta^\epsilon (t),h^\epsilon (t),\xi^{-,\epsilon}_k(t)),\\
q_0(t)&=(\xi^{+,0}_k(t),J^0(t),\theta^0(t),h^0(t),\xi^{-,0}_k(t)).\end{align*}
Let $T_1(>T)$ be a time such that
\begin{equation} \| q_\epsilon (t)-q_0(t)\| _n\leq \tilde{C}_2 \epsilon 
|\ln \epsilon |^2,\end{equation}
for all $t\in [T,T_1]$, where $\tilde{C}_2 =\tilde{C}_2(T)$ is 
independent of $\epsilon$. 
 From \eqref{sme1}, such a $T_1$ exists. The proof will be completed
through a continuation argument. For $t\in [T,T_1]$,
\begin{equation}\begin{split}\label{smcd}&\sum_{k=1,2} 
(|\xi^{+,0}_k(t)|+|\xi^{-,0}_k(t)|)+\|h^0(t)\|_n\leq 
C_3re^{-\frac{1}{2}\mu
(t-T)},\quad \ |J^0(t)|\leq C_4\sqrt{\epsilon},\\
& |J^\epsilon (t)|\leq |J^0(t)|+|J^\epsilon(t)-J^0(t)| \leq 
|J^0(t)|+\tilde{C}_2\epsilon |\ln
\epsilon |^2\leq C_5\sqrt{\epsilon},\\
&\sum_{k=1,2} 
(|\xi^{+,\epsilon}_k(t)|+|\xi^{-,\epsilon}_k(t)|)+\|h^\epsilon 
(t)\|_n\leq C_3re^{-\frac{1}{2}\mu (t-T)} +\tilde{C}_2\epsilon
|\ln \epsilon |^2,\end{split}\end{equation}
where $r$ is small. Since actually $q_\epsilon (t)$, $q_0(t)\in H^n$ 
for any fixed $n\geq 1$, by Theorem~\ref{CSM},
\begin{equation} |\xi^{+,\epsilon}_k(t)-\xi^{+,0}_k(t)|\leq C_6\| 
v_\epsilon (t)-v_0(t)\|_n+C_7\epsilon ,\end{equation}
whenever $v_\epsilon (t)$, $v_0(t)\in E_{n+4}(r)$, where $v_\epsilon 
(T)=v_\epsilon $ and $v_0(T)=v_0$ are defined in \eqref{sme}. Thus
we only need to estimate $\| v_\epsilon (t)-v_0(t)\|_n$. From 
\eqref{cspe2}, we have for $t\in [T, T_1]$ that
\begin{equation}\label{smeq1} 
v(t)=e^{A(t-T)}v(T)+\int^t_Te^{A(t-\tau)}\tilde{F}(\tau 
)d\tau.\end{equation}
Let $\Delta v(t)=v_\epsilon (t)-v_0(t)$. Then
\begin{equation}\begin{split} \Delta 
v(t)&=[e^{A(t-T)}-e^{A\mid_{\epsilon 
=0}(t-T)}]v_0(T)+e^{A(t-T)}\Delta v(T)\\
&\quad + \int^t_Te^{A(t-\tau )}[\tilde{F}(\tau )-\tilde{F}(\tau 
)|_{\epsilon =0}]d\tau \\
&\quad + \int^t_T[e^{A(t-\tau )}-e^{A\mid_{\epsilon =0}(t-\tau 
)}]\ \tilde{F}(\tau )|_{\epsilon =0}d\tau.\end{split}\end{equation}
By the condition \eqref{smcd}, we have for $t\in [T, T_1]$ that
\begin{equation}\| 
\tilde{F}(t)-\tilde{F}(t)|_{\epsilon=0}\|_n\leq 
[C_8\sqrt{\epsilon}+C_9 re^{-\frac{1}{2}\mu (t-T)} ]\epsilon |\ln
\epsilon |^2.\end{equation}
Then
\begin{equation} \| \Delta v(t)\|_n\leq C_{10}\epsilon 
(t-T)+C_{11}r\epsilon |\ln \epsilon 
|^2+C_{12}\sqrt{\epsilon}(t-T)^2\epsilon
|\ln \epsilon |^2.\end{equation}
Thus by the continuation argument, for $t\in [T, T+\frac{1}{\mu}|\ln 
\epsilon |]$, there is a constant $\hat{C}_1=\hat{C}_1(T)$,
\begin{equation} \| \Delta v(t)\|_n\leq \hat{C}_1\epsilon |\ln 
\epsilon |^2.\end{equation}
\end{proof}
By Lemma~\ref{lesm} and estimate \eqref{sarg1},
\begin{equation}\label{sarg2} \text{distance} \left\{ q_\epsilon 
( T+\frac{1}{\mu}|\ln \epsilon |), \Pi\right\}<\tilde{C}\epsilon |\ln \epsilon |^2.\end{equation}
Recall the fish-like singular level set given by $\mathcal{H}$ 
\eqref{fham}, the width of the fish is of order
$\mathcal{O}(\sqrt{\epsilon })$, and the length of the fish is of 
order $\mathcal{O}(1)$. Notice also that $q_0(t)$ has a phase shift
\begin{equation} \theta^0_1=\theta^0( T+\frac{1}{\mu}|\ln 
\epsilon |)-\theta^0(0).\end{equation}
For fixed $\beta$, changing $\alpha$ can induce $\mathcal{O}(1)$ 
change in the length of the fish, $\mathcal{O}(\sqrt{\epsilon})$ 
change
in $\theta^0_1$, and $\mathcal{O}(1)$ change in $\theta^0(0)$. See 
Figure~\ref{supsm} for an illustration.
The leading order signed distance from $q_\epsilon ( 
T+\frac{1}{\mu}|\ln \epsilon |)$ to $W^s_n(Q_\epsilon )$ can be 
defined
as
\begin{equation}\begin{split} \tilde{d}
&=\mathcal{H}(j_0,\theta^0(0))-\mathcal{H}(j_0,\theta^0(0)+\theta^0_1)\\
&=2\omega \left [\alpha \omega \theta^0_1+\beta [\sin \theta^0(0)-\sin 
(\theta^0(0)+\theta^0_1)]\right ],\end{split}\label{dtd}\end{equation}
where $\mathcal{H}$ is given in \eqref{fham}. The common zero of 
$M_k$ \eqref{Meln1} and $\tilde{d}$ and the implicit function 
theorem
imply the existence of a homoclinic orbit asymptotic to $Q_\epsilon$. 
Much detailed arguments have been given to the signed distance
$\tilde{d}$ and the second measurement in \cite{LMSW96} and \cite{LM97}.

\begin{figure}
\vspace{1.5in}
\caption{The second measurement.}
\label{supsm}
\end{figure}

\subsection{Evaluation of Melnikov Integrals and Second Distance}

It turns out that to the leading order, one can evaluate $M_k$
(\ref{Meln1}) at $q_0(t)$ where $a =\om$. Our goal in this subsection 
is to find the common zero of $M_k$ (\ref{Meln1}) and $\tilde{d}$ 
(\ref{dtd}). \\

\nid
{\bf One Pair of Complex Double Points Case}

$M_1 = 0$ and $\tilde{d}=0$ lead to 
\begin{equation}
M_1 = M^{(1)} + \al M^{(2)} + \be \cos \ga M^{(3)} = 0\ , 
\label{ecml1}
\end{equation}
\begin{equation}
\be \cos \ga = \frac {\al \om (\Dl \ga )} {2 \sin \frac {\Dl \ga }{2}} \ ,
\label{ecml2}
\end{equation}
where $\Dl \ga = -4 \vth_0$, $M^{(j)} = M^{(j)}(\om)$, ($j=1,2,3$), 
and 
\[
M^{(1)} = \om^2 \int^{+\infty}_{-\infty}\int^{2\pi}_{0} 
(|u_1|^2+|u_2|^2)^{-2}[\bar{u}_1^2 \pa_x^2 P -  
\bar{u}_2^2 \pa_x^2 \bar{P}] dx dt \ ,
\]
\[
M^{(2)} = \om^2 \int^{+\infty}_{-\infty}\int^{2\pi}_{0} 
(|u_1|^2+|u_2|^2)^{-2}[\bar{u}_2^2 \bar{P}-  
\bar{u}_1^2 P] dx dt \ ,
\]
\[
M^{(3)} = \om \int^{+\infty}_{-\infty}\int^{2\pi}_{0} 
(|u_1|^2+|u_2|^2)^{-2}[\bar{u}_1^2 -  
\bar{u}_2^2 ] dx dt \ ,
\]
and $P$ is given by $Q=q_c P$, and $Q$ is given in (\ref{se}). 
Equations (\ref{ecml1}) and (\ref{ecml2}) define a codimension-one
surface in the space of ($\al,\be,\om$), given by
\[
\al = \frac {1}{\k(\om)} \ ,
\]
where
\[
\k(\om) = - [2 M^{(2)} \sin \frac {\Dl \ga }{2} + M^{(3)}  
\om (\Dl \ga)][2 M^{(1)} \sin \frac {\Dl \ga }{2}]^{-1} \ ,
\]
and its graph is plotted in Figure \ref{opcfig}. \\

\nid
{\bf Two Pairs of Complex Double Points Case}

$M_j = 0$ ($j=1,2$) and $\tilde{d}=0$ lead to 
\begin{equation}
M_j = M_j^{(1)} + \al M_j^{(2)} + \be \cos \ga M_j^{(3)} 
+ \be \sin \ga M_j^{(4)}= 0\ ,  \ \ (j=1,2) 
\label{ecml3}
\end{equation}
\begin{equation}
\be \cos \ga = \frac {\al \om \widetilde{\Dl \ga }} {2 \sin \frac 
{\widetilde{\Dl \ga} }{2}} \ ,
\label{ecml4}
\end{equation}
where $\widetilde{\Dl \ga } = -4 (\vth_0 +\hvth_0)$, $M_j^{(l)} = 
M_j^{(l)}(\om, \Dl \rho)$, ($j=1,2, \ l=1,2,3,4$), $\Dl \rho =
2 \hat{\sg} \sg^{-1} \rho - \hat{\rho}$, and 
\begin{eqnarray*}
M_1^{(1)}&=& \om^2 \int^{+\infty}_{-\infty}\int^{2\pi}_{0} 
[S_2^2 \pa_x^2 \tilde{P} - S_1^2 \pa_x^2 \overline{\tilde{P}}]dxdt\ ,
\\
M_1^{(2)}&=& \om^2 \int^{+\infty}_{-\infty}\int^{2\pi}_{0} 
[S_1^2 \overline{\tilde{P}} - S_2^2 \tilde{P}]dxdt\ ,
\\
M_1^{(3)}&=& \om \int^{+\infty}_{-\infty}\int^{2\pi}_{0} 
[S_2^2 - S_1^2]dxdt\ ,
\\
M_1^{(4)}&=& i\om \int^{+\infty}_{-\infty}\int^{2\pi}_{0} 
[S_2^2 + S_1^2]dxdt\ ,
\end{eqnarray*}
and $\tilde{P}$ is given by $\tilde{Q}= q_c \tilde{P}$, and $\tilde{Q}$
is given in (\ref{dne}) and (\ref{dec}). $M_2^{(l)}$ can be obtained 
from $M_1^{(l)}$ ($l = 1,2,3,4$) by replacing $S_m$ by $\hat{S}_m$ 
($m=1,2$). 
\begin{figure}
\vspace{1.5in}
\caption{The curve of $\k (\om )$.}
\label{opcfig}
\end{figure}
\begin{figure}
\vspace{1.5in}
\caption{The surface of $\tilde{\chi}(\om, \Dl \rho)$.}
\label{dpcfig}
\end{figure}
\nid
Equations (\ref{ecml3}) and (\ref{ecml4}) define a codimension-one
surface in the space of ($\al,\be,\om$), given by
\begin{eqnarray*}
\al &=& \frac {1}{\tilde{\chi}(\om, \Dl \rho)}\ , \\
\be &=& \be (\om, \Dl \rho) = \bigg [ (\al \om \widetilde{\Dl \ga})^2 
(2 \sin \frac {\widetilde{\Dl \ga}}{2})^{-2} \\
& & + (M_1^{(4)})^{-2} 
[M_1^{(1)} + \al (M_1^{(2)} + M_1^{(3)}\om \widetilde{\Dl \ga} 
(\sin \frac {\widetilde{\Dl \ga}}{2})^{-1})]^2 \bigg ]^{1/2}\ ,
\end{eqnarray*}
where
\begin{eqnarray*}
\tilde{\chi}(\om, \Dl \rho)&=& (M_2^{(1)} M_1^{(4)}- M_1^{(1)}
M_2^{(4)})^{-1}\bigg [ (M_1^{(2)} M_2^{(4)}- M_2^{(2)}
M_1^{(4)}) \\
& & +\om \widetilde{\Dl \ga }(2 \sin \frac 
{\widetilde{\Dl \ga}}{2})^{-1}
(M_1^{(3)} M_2^{(4)}- M_2^{(3)}M_1^{(4)}) \bigg ]\ ,
\end{eqnarray*}
and its graph is plotted in Figure \ref{dpcfig}.

\subsection{Statement of the main Theorem}

\begin{theorem}[Main Theorem] There exists a $\e_0 > 0$, such that 
for any $\e \in (0, \e_0)$, there exists 
a codimension 1 surface in the space of $(\alpha,\beta, \om) \in 
\mathbb{R}^+\times  \mathbb{R}^+\times \mathbb{R}^+$ where 
$\om \in (\frac{1}{2}, \frac{3}{2})/S$, $S$ is a finite subset, and 
$\al \om < \be$. For any $(\alpha ,\beta, \omega)$ on the codimension-one
surface, the perturbed nonlinear Schr\"odinger equation 
\eqref{pnls} possesses a homoclinic orbit asymptotic to the saddle
$Q_\epsilon$ \eqref{Qec}. The codimension 1 surface has the 
approximate representation given in the 
subsection on Evaluation of Melnikov Integrals and Second 
Measurement.\end{theorem}
\begin{proof}
From the explicit computation in last subsection and the implicit 
function theorem, $d_k$ and $\tilde{d}$ are zero for the parameter 
values specified in the theorem .
\end{proof}